\documentclass[secthm,seceqn,amsthm,ussrhead,reqno, 12pt]{amsart}
\usepackage[utf8]{inputenc}
\usepackage[english]{babel}
\usepackage[symbol]{footmisc}
\usepackage{amssymb,amsmath,amsthm,amsfonts,xcolor,enumerate,hyperref,comment,longtable,cleveref}

\usepackage{times}
\usepackage{cite}
\usepackage{pdflscape}
\usepackage{ulem}
\usepackage[mathcal]{euscript}
\usepackage{tikz}
\usepackage{hyperref}
\usepackage{cancel}
\usepackage{stmaryrd}

\newtheorem{thrm}{Theorem}

\newtheorem{corollary}[thrm]{Corollary}

\newtheorem{proposition}[thrm]{Proposition}
\newtheorem{remark}[thrm]{Remark}

\newtheorem{defn}[thrm]{Definition}
\newtheorem{definition}[thrm]{Definition}

\crefrangeformat{equation}{#3(#1)#4--#5(#2)#6}

\crefname{thrm}{thrm}{Theorems}
\crefname{lem}{Lemma}{Lemmas}
\crefname{cor}{Corollary}{Corollaries}
\crefname{prop}{Proposition}{Propositions}
\crefname{proposition}{Proposition}{Propositions}
\crefname{defn}{Definition}{Definitions}
\crefname{exm}{Example}{Examples}
\crefname{rem}{Remark}{Remarks}
\crefname{section}{Section}{Sections}
\crefname{equation}{\unskip}{\unskip}
\crefname{enumi}{\unskip}{\unskip}

\begin{document}

\noindent{\Large
Local and 2-local $\frac{1}{2}$-derivations on
 finite-dimensional Lie algebras}

	\bigskip
	
	 \bigskip

\begin{center}	
{\bf Abror Khudoyberdiyev\footnote{V.I.Romanovskiy Institute of Mathematics Academy of Science of Uzbekistan; National University of Uzbekistan; \ khabror@mail.ru}},	
{\bf
Bakhtiyor Yusupov\footnote{ V.I.Romanovskiy Institute of Mathematics Academy of Science of Uzbekistan; \
Department of Physics and Mathematics, Urgench State University;
 \
baxtiyor\_yusupov\_93@mail.ru}}
\end{center}

	\bigskip

\noindent {\bf Abstract.}
{\it In this work, we introduce the notion of local and $2$-local $\delta$-derivations and describe local and $2$-local $\frac{1}{2}$-derivation of finite-dimensional solvable Lie algebras with filiform, Heisenberg, and abelian nilradicals. Moreover, we describe the local $\frac{1}{2}$-derivation of oscillator Lie algebras, Schr{\"o}dinger algebras, and Lie algebra with a three-dimensional simple part, whose radical is an irreducible module. We prove that an algebra with only trivial $\frac{1}{2}$-derivation does not admit local and $2$-local $\frac{1}{2}$-derivation, which is not $\frac{1}{2}$-derivation.
}

\bigskip



 \bigskip

\bigskip
\section{Introduction}

The notion of $\delta$-derivations was initiated by V.Filippov for Lie algebras in \cite{fil1, fil2}. The space
of $\delta$-derivations includes usual derivations ($\delta=1$), anti-derivations ($\delta=-1$) and elements from the centroid. In \cite{fil2} it was proved that prime Lie algebras, as a rule, do not have nonzero $\delta$-derivations (in case of $\delta \neq 1, -1, 0, \frac 1 2$),
and all $\frac 1 2$-derivations of an arbitrary prime Lie algebra $A$ over the field $\mathbb{F}$ ($\frac 1 6 \in \mathbb{F}$) with a non-degenerate symmetric invariant bilinear form were described. It was proved that if $A$ is a central simple Lie algebra over a field
of characteristic $p \neq 2, 3$ with a non-degenerate symmetric invariant bilinear form, then any $\frac 1 2$-derivation $D$ has the form $D(x) = \lambda x$ for some $\lambda \in \mathbb{F}.$

In \cite{fil3}, $\delta$-derivations were investigated for prime alternative and non-Lie Mal’tsev algebras, and it was proved that alternative and non-Lie Mal’tsev algebras with certain restrictions on the ring of operators $F$ have no non-trivial $\delta$-derivation. A description of  $\delta$-derivations of classical Lie superalgebras was given in \cite{Kay2}. $\delta$-derivations of finite-dimensional semi-simple Jordan algebras over the field of characteristic different from $2$ and $\delta$-
superderivations of finite-dimensional simple Lie and Jordan superalgebras were investigated in \cite{Kay3}. P.Zusmanovich in \cite{Zus} described $\delta$-(super)derivations of prime Lie superalgebras, namely, he proved that a prime Lie superalgebra has no non-trivial $\delta$-(super)derivations for $\delta\neq  1, -1, 0, \frac 1 2$.
%Note that $\delta$-derivations is a particular case of generalized derivations (see, \cite{Leger}).

Nowadays, local and $2$-local operators have become popular for some non-associative algebras such as the Lie, Jordan, and Leibniz algebras. The notions of local derivations were introduced in 1990 by R.Kadison \cite{Kadison} and D.Larson, A.Sourour \cite{Larson}. Later in 1997, P.\v{S}emrl introduced the notions of  $2$-local derivations and $2$-local automorphisms on algebras \cite{Sem}. 

Investigation of local derivations on Lie algebras was initiated in  \cite{Ayupov7} by Sh.Ayupov and K.Kudaybergenov. They proved that every local derivation on semi-simple Lie algebras is a derivation and gave examples of nilpotent finite-dimensional Lie algebras with local derivations that are not derivations. In \cite{AyuKhud}, local derivations of solvable Lie algebras are investigated, and it is shown that in the class of solvable Lie algebras there exist algebras that admit local derivations that are not derivations and also algebras for which every local derivation is a derivation. Several authors investigated local derivations for the finite or infinite dimensional Lie and Leibniz algebras \cite{KOK22, YC20, CZZ, Yao, AKYu1, AEK1, AY1, AyuKhudYus, KayKhudYul, AKO}. It was proved that 
 all local derivations of the following algebras are derivations:
  Borel subalgebras of finite-dimensional simple Lie algebras; infinite-dimensional Witt algebras over an algebraically closed field of characteristic zero; solvable Lie algebras of maximal rank; Cayley algebras;  
  locally finite split simple Lie algebras over a field of characteristic zero; the Schr\"{o}dinger algebras; conformal Galilei algebras.

Several papers have been devoted to similar notions and corresponding problems for $2$-local derivations and automorphisms of Lie algebras \cite{AKR, ChenWang, AKY, XMP, AKYu1, QX, YK, ChenXe }. Namely, in \cite{AKR} it is proved that every $2$-local derivation on the semi-simple Lie algebras is a derivation and that
each finite-dimensional nilpotent Lie algebra, with dimension
larger than two admits $2$-local derivation, which is not a
derivation.  Let us present a list of finite or infinite-dimensional  Lie algebras for which all 2-local derivations are derivations: finite-dimensional semi-simple Lie algebras over an algebraically closed field of characteristic zero; infinite-dimensional Witt  algebras over an algebraically closed field of characteristic zero;
		 locally finite split simple Lie algebras over a field of characteristic zero;
        Virasoro algebras; Virasoro-like algebra;
        the Schrödinger-Virasoro algebra;
         Jacobson-Witt algebras; planar Galilean conformal algebras.

In \cite{BKL} by Ferreira, Kaygorodov, Lopatkin,
	a relation between $\frac{1}{2}$-derivations of Lie algebras and 
	transposed Poisson algebras has been established. This relation was used to describe all transposed Poisson structures of several classes of Lie algebras. Descriptions of $\frac{1}{2}$-derivations	are given on  Witt and Virasoro algebras in  \cite{BKL};
	on   twisted Heisenberg-Virasoro,  Schr\"odinger-Virasoro,  and  
	extended Schr\"odinger-Virasoro algebras in \cite{yh21};
	on Schr\"odinger algebra in $(n+1)$-dimensional space-time in \cite{ytk};
on solvable Lie algebra with filiform, Heisenberg, and abelian nilradicals in \cite{KayKhud1};
	on Witt type algebras in \cite{kk23};
	on generalized Witt algebras and  Block Lie algebras in \cite{kkg23}.

In this work, we introduce the notion of local and $2$-local $\delta$-derivations and prove that an algebra with only trivial $\frac{1}{2}$-derivation does not admit local and $2$-local $\frac{1}{2}$-derivation, which is not $\frac{1}{2}$-derivation.
We describe the local and $2$-local $\frac{1}{2}$-derivation of finite-dimensional solvable Lie algebras with filiform, Heisenberg, and abelian nilradicals. Moreover, local and $2$-local $\frac{1}{2}$-derivation of oscillator Lie algebras, Schr{\"o}dinger algebras, and Lie algebra with a three-dimensional simple part, whose radical is an irreducible module are described.
   
\section{Preliminaries}\label{prem}

All the algebras below will be over the complex field, and all the linear maps will be $\mathbb C$-linear unless otherwise stated.
 
\begin{defn}\label{12der}
		Let $({\mathfrak L}, [-,-])$ be an algebra with a multiplication $[-,-].$ A linear map $D$ is called a $\delta$-derivation if it satisfies
\begin{center}
$D[x,y]= \delta \big([D(x),y]+ [x, D(y)] \big).$
\end{center}
	\end{defn}

Note that $1$-derivation is a usual derivation and $(-1)$-derivation is called an anti-derivation.
If $D_1$ and $D_2$ are $\delta_1$ and $\delta_2$-derivations, respectively, then their commutator $[D_1, D_2] = D_1D_2 - D_2D_1$ is a $\delta_1\delta_2$- derivation.   The set of all $\delta$-derivations, for the fixed $\delta$, we denote by $Der_{\delta}(\mathfrak L).$

%A vector space ${\mathfrak L}$ over a field $\mathbb{F},$ with an operation ${\mathfrak L}\times {\mathfrak L}\rightarrow {\mathfrak L},$ denoted $(x,y)\mapsto[x,y]$ and called the \textbf{bracket} or \textbf{commutator} of $x$ and $y,$ is called a \textbf{Lie algebras} over $\mathbb{F}$ if the following axioms are satisfied:

%\begin{itemize}
%  \item The bracket operation is bilinear.
%  \item $[x,x] =0 ,\ \ \ \  \forall x\in{\mathfrak L}.$
%  \item $[[x,y],z]+[[y,z],x]+[[z,x],y]=0, \forall x,y,z\in{\mathfrak L}.$
%\end{itemize}

%Observe that $\frac{1}{2}$-derivations are a particular case of $\delta$-derivations .
%
%It is easy to see from Definition \ref{12der} that $[\LL,\LL]$ and $\Ann(\LL)$ are invariant under any $\frac 12$-derivation of $\LL$.

For the fixed $\delta$, we give the notions of local and $2$-local $\delta$-derivation.

\begin{definition}
A linear operator $\Delta$ is called a local $\delta$-derivation, if for any $x \in \mathfrak L,$ there exists a $\delta$-derivation $D_x: \mathfrak L \rightarrow \mathfrak L$ (depending on $x$) such that $\Delta(x) = D_x(x).$ The set of all local $\delta$-derivations on $\mathfrak L$ we denote by $\mathrm{Loc}\mathrm{Der}_{\delta}(\mathfrak L).$
\end{definition}

\begin{definition}
A map $\nabla :  \mathfrak L \rightarrow \mathfrak L$ (not
necessary  linear) is called a $2$-local $\delta$-derivation, if for any $x,y\in \mathfrak L$, there exists a $\delta$-derivation $D_{x,y}\in \mathrm{Der}_{\delta}
(\mathfrak L)$ such that
{\small\[
\nabla(x)=D_{x,y}(x), \quad \nabla(y)=D_{x,y}(y).
\]}
\end{definition}

It should be noted that $2$-local $\delta$-derivation is not necessary linear, but for any $x \in \mathfrak L$ and for any scalar $\lambda,$ we have that
$$\nabla(\lambda x) = D_{x, \lambda x}(\lambda x) = \lambda  D_{x, \lambda x}(x) = 
\lambda \nabla(x).$$

In this work, we focus on investigating local and $2$-local $\frac 1 2$-derivations. 
Note that the main example of $\frac 1 2$-derivations is the multiplication by an element from the ground field, i.e., $D(x) = \lambda x$ for all $x \in \mathfrak L$.
Such kind of $\frac 1 2$-derivations are called trivial $\frac 1 2$-derivations.
In \cite{BKL} it was proved that 
if a Lie algebra $\mathfrak L$ does not admit non-trivial $\frac{1}{2}$-derivations, then all transposed Poisson algebra structures on $\mathfrak L$ are trivial.

In the following theorem, we show that if an algebra $\mathfrak L$ does not admit non-trivial $\frac{1}{2}$-derivations, then any local and $2$-local $\frac{1}{2}$-derivation of the algebra $\mathfrak L$ is a $\frac{1}{2}$-derivation. 
\begin{thrm}\label{thm1}
Let $\mathfrak{L}$ be an algebra, whose all $\frac{1}{2}$-derivation are trivial. Then any local and $2$-local $\frac{1}{2}$-derivation of $\mathfrak{L}$ is a trivial $\frac{1}{2}$-derivation.
\end{thrm}

\begin{proof} We can assume $dim{\mathfrak{L}} \geq 2$, since the theorem is evident for the one-dimensional case.

Let $\Delta$ be a local $\frac{1}{2}$-derivation on $\mathfrak L$, then for any linear independent elements $x$ and $y$ of $\mathfrak L$, there exist $\frac{1}{2}$-derivations $D_{x}$ and $D_{y}$ such that 
$$\Delta(x) = D_x(x) = \alpha_x x, \quad \Delta(y) = D_y(y) = \alpha_y y,$$
which implies $\Delta(x+y) = \alpha_x x+ \alpha_y y.$

On the other hand, for the element $x+y$, we have a $\frac{1}{2}$-derivation $D_{x+y}$ such that $$\Delta(x+y) = D_{x+y}(x+y) = \alpha_{x+y}(x+y)= \alpha_{x+y}x+\alpha_{x+y}y.$$

Comparing  the coefficients of $x$ and $y$, we obtain that $\alpha_x =\alpha_y.$
Due to the arbitrariness of $x$ and $y$, we get that $\Delta(x) = \alpha x$ for any element $x \in \mathfrak L.$ Therefore, $\Delta$ is a
$\frac{1}{2}$-derivation.

Now, let $\nabla$ be a $2$-local $\frac{1}{2}$-derivation on $\mathfrak L.$ Then for the fixed element $x$ and for arbitrary elements $y$ and $z$, we have that there exist  $\frac{1}{2}$-derivations $D_{x,y}$ and $D_{x,z}$, such that
$$\nabla(x)=D_{x,y}(x) = \alpha_{x,y} x,\quad \nabla(y)=D_{x,y}(y) = \alpha_{x,y} y,$$
$$\nabla(x)=D_{x,z}(x) = \alpha_{x,z} x,\quad \nabla(z)=D_{x,z}(z) = \alpha_{x,z} z.$$

Hence, we obtain $\alpha_{x,y} = \alpha_{x,z}.$ Due to the arbitrariness of $y$ and $z$, we get that $\nabla(y) = \alpha y$ for any element $y\in \mathfrak L.$ Therefore, $\nabla$ is a $\frac{1}{2}$-derivation.
\end{proof}

It should be noted that the definition of $\frac 1 n$-derivations of $n$-ary algebras as a particular case of $\delta$-derivations of $n$-ary algebras was given in \cite{Kay-n-ary}.

\begin{definition} Let $(\mathfrak{L}, [ , \dots , ])$ be an $n$-ary algebra with the multiplication $[ , \dots ,]$ and $\varphi$
be a linear map. Then $\varphi$ is a $\frac 1 n$-derivation if it satisfies
$$\varphi([x_1,x_2,\dots,x_n]) = \frac 1 n \sum\limits_{i=1}^n[x_1, \dots, \varphi(x_i), \dots, x_n].$$
\end{definition}

The main example of $\frac 1 n$-derivations is the multiplication of an element from the basic field, and such $\frac 1 n$-derivations are called trivial $\frac 1 n$-derivations.

Theorem \ref{thm1} is also true for the $n$-ary case, that is, if a $n$-ary algebra has only a trivial $\frac{1}{n}$-derivation, then any local and $2$-local $\frac{1}{n}$-derivation is a $\frac{1}{n}$-derivation. 
Moreover, the results of Theorem \ref{thm1} can be easily extended for the algebras, which are direct sums of the vector spaces, and all $\frac{1}{n}$-derivation invariant and trivial on these subspaces. 
For example, every complex finite-dimensional semi-simple Lie algebra is a direct sum of simple algebras, and any $\frac{1}{2}$-derivation will be invariant on all simple algebras from this direct sum. From this observation and the results of the works \cite{fil2, fil3, Kay3, Kay-n-ary, KayOpa, kk23, kkg23, klv22, BKL} and Theorem \ref{thm1}, we have the following corollary:

\begin{corollary} Any local and $2$-local $\frac{1}{2}$-derivation of the following algebras is a $\frac{1}{2}$-derivation:
\begin{itemize}
\item finite-dimensional semi-simple Lie algebras;

\item finite-dimensional semi-simple Jordan algebras;

\item  finite-dimensional semi-simple Mal'tsev algebras; 

\item finite-dimensional semi-simple structurable algebras;

\item finite-dimensional semi-simple alternative algebras;

\item finite-dimensional semi-simple  $n$-Lie algebras;

\item The Lie algebra $\mathcal{W}(a, b)$ for $b\neq -1,$ with basis $\{L_i , I_j\}_{i,j \in \mathbb{Z}}$ and multiplication
$$[L_m, L_n] = (m - n)L_{m+n}, [L_m, I_n] = -(n + a + bm)I_{m+n}.$$

\item The Virasoro algebra \textbf{Vir}  with basis $\{C, L_i\}_{i\in \mathbb{Z}}$ and multiplication
$$[L_m, L_n] = (m - n)L_{m+n} + \frac{m^3-m}{12}\delta_{m+n, 0}C.$$

\item The Block Lie algebra $B(q)$ for $q \neq \mathbb{Z}$ with basis $\{L_{m,i}\}_{m, i\in \mathbb{Z}}$ and multiplication
$$[L_{m,i}, L_{n,j}] = (n(i+q) - m(j+q))L_{m+n,i+j}.$$

\item Galilean algebras. 

%\item conformal centrally extended Galilei algebra $\widetilde{\mathfrak{g}}^{(\ell)}.$
%$$\begin{array}{lll}[h, e] = 2e, & [h, f ] = -2 f, & [e, f ] = h, \\[1mm][h, p_k] = 2(\ell - k)p_k, & [e, p_k] = kp_{k-1}, & [f, p_k] = (2\ell - k)p_{k+1}, \\[1mm] \multicolumn{3}{l} {[p_k, p_{2\ell-k}] = (-1)^{k+\ell+\frac 12} k!(2\ell - k)!z }  \end{array}$$
%where $0 < \ell \in \mathbb{N} - \frac 12,$ $0 \leq k \leq 2 \ell.$

\end{itemize}

\end{corollary}

%It is proven that any finite-dimensional solvable Lie algebra over an algebraically closed field of zero characteristic admits non-trivial $\frac12$-derivations \cite{klv22}.

\section{Local and 2-local $\frac{1}{2}$-derivation on  finite-dimensional solvable Lie algebras}\label{solv}

For a Lie algebra ${\mathfrak L}$, consider the following lower central and
derived sequences:
$$
{\mathfrak L}^1={\mathfrak L},\quad {\mathfrak L}^{k+1}=[{\mathfrak L}^k,{\mathfrak L}^1], \quad k \geq 1,
$$
$${\mathfrak L}^{[1]} = {\mathfrak L}, \quad {\mathfrak L}^{[s+1]} = [{\mathfrak L}^{[s]}, {\mathfrak L}^{[s]}], \quad s \geq 1.$$

\begin{definition} A Lie algebra ${\mathfrak L}$ is called
nilpotent (respectively, solvable), if there exists  $p\in\mathbb N$ $(q\in
\mathbb N)$ such that ${\mathfrak L}^p=0$ (respectively, ${\mathfrak L}^{[q]}=0$). The minimal number $p$ (respectively, $q$) with such
property is said to be the index of nilpotency (respectively, of solvability) of the algebra ${\mathfrak L}$.
\end{definition}
Note that any Lie algebra ${\mathfrak L}$ contains a unique maximal solvable (resp. nilpotent) ideal, called the radical (resp. nilradical) of the algebra.

Note that, any finite-dimensional solvable Lie algebra over an algebraically closed field of zero characteristic admits non-trivial $\frac12$-derivations \cite{klv22}. 
In this section, we describe local and 2-local $\frac{1}{2}$-derivation of oscillator Lie algebras, and solvable Lie algebras with naturally graded filiform, Heisenberg, and abelian nilradicals.

\subsection{Local and 2-local $\frac{1}{2}$-derivation on solvable Lie algebras with naturally graded filiform nilradical}

It is well known that there are two types of naturally graded filiform Lie algebras. The second type will appear only in the case when the dimension of the algebra is even.
 Any naturally graded filiform Lie algebra is isomorphic to one of the following non-isomorphic algebras \cite{Ver}:
\begin{longtable}{llllll}
$n_{n,1}$ &:& $[e_i, e_1]=-[e_1, e_i]=e_{i+1},$ & $2\leq i \leq n-1.$\\ \hline
$Q_{2n}$&:& $[e_i, e_1] =  -[e_1, e_i]=e_{i+1}$,& $2\leq i \leq 2n-2,$\\
&& $[e_i, e_{2n+1-i}]  =  -[e_{2n+1-i}, e_i]=(-1)^i e_{2n}$,& $2\leq i\leq n.$
   \end{longtable}

All solvable Lie algebras whose nilradical is the naturally graded filiform Lie algebras $n_{n,1}$ and $Q_{2n}$ are classified in \cite{AnCaGa3}, \cite{sw05}. 
Here we give the list of solvable Lie algebras with nilradical $n_{n,1}$:
\begin{longtable}{llllll}
\hline
$\mathfrak{s}^{1}_{n,1}(\beta)$ & : &
$[e_i, e_1] = e_{i+1}$, &   $2 \leq i \leq n-1,$ \\
& & $[e_i, x]  =(i-2+\beta) e_i$,  & $2 \leq i \leq n$, & $[e_1,x]=e_1$. \\ \hline
$\mathfrak{s}^{2}_{n,1}$ & : & 
$[e_i, e_1] = e_{i+1}$, &  $2 \leq i \leq n-1$, \\
 & & $[e_i, x]  = e_i$,  &  $2 \leq i \leq n.$ \\ \hline
$\mathfrak{s}^{3}_{n,1}$  & : & 
$[e_i, e_1] = e_{i+1}$, &  $2 \leq i \leq n-1,$ \\
& & $[e_i, x]  = (i-1)e_i$,  & $2 \leq i \leq n,$ & $[e_1,x]=e_1 +e_2$.
 \\ \hline
 $\mathfrak{s}^{4}_{n,1}(\alpha_3, \dots, \alpha_{n-1})$ & : &
 $[e_i,e_1]=e_{i+1}$, & $2 \leq i \leq n-1$,\\
 & & $[e_i,x]=e_i+\sum\limits_{l=i+2}^{n} \alpha_{l+1-i} e_l$,  & $2\leq i\leq n.$ \\ \hline
$\mathfrak{s}_{n,2}$ & : &  $[e_i, e_1] = e_{i+1}$, &  $2 \leq i \leq n-1$, \\
& & $[e_i, x_1]  = (i-2)e_i$,  & $3 \leq i \leq n,$ & $[e_1, x_1] = e_1,$\\
 && $[e_i, x_2] = e_i$, & $2 \leq i \leq n.$ \\ \hline
   \end{longtable}

The list of solvable Lie algebras with nilradical $Q_{2n}$:
\begin{longtable}{lllll}
\hline
$\tau_{2n,1}^1(\alpha)$ & : & 
$[e_i,e_1] =e_{i+1}$, & $2\leq i \leq 2n-2,$ \\
&&$[e_i,e_{2n+1-i}] =(-1)^i e_{2n},$ & $2\leq i \leq n,$ &
 $[e_1,x]=e_1$,\\
 && $[e_i,x]=(i-2+\alpha) e_i$, & $2\leq i\leq 2n-1$, &
$[e_{2n},x] =(2n-3+2\alpha)e_{2n}.$ \\
\hline
$\tau_{2n,1}^2$ & : & 
$[e_i,e_1] =e_{i+1}$, & $2\leq i \leq 2n-2,$\\
& & $[e_i,e_{2n+1-i}] =(-1)^i e_{2n}$,& $2\leq i \leq n$, &
$[e_1,x]=e_1+ e_{2n}$,\\
 & & $[e_i,x]=(i-n) e_i$, & $2\leq i\leq 2n-1$, &
$[e_{2n},x] = e_{2n}$.\\
\hline
\multicolumn{5}{l}{$\begin{array}{llll} 
\tau_{2n,1}^3(\alpha_4, \alpha_6, \dots \alpha_{2n-2})& : &
[e_i,e_1]=e_{i+1}, & 2\leq i \leq 2n-2,\\
&& [e_i,e_{2n+1-i}]=(-1)^i e_{2n},& 2\leq i \leq n,\\
 && [e_{i+2},x] =e_{i+2} + \sum\limits_{k=2}^{\lfloor \frac {2n-3-i} {2} \rfloor}\alpha_{2k} e_{2k+1+i}, & 0\leq i \leq 2n-3, \\ &&
  [e_{2n},x] = 2 e_{2n}.\end{array}$}\\
\hline
$\tau_{2n,2}$ & : &
$[e_i,e_1]=e_{i+1}$, & $2\leq i \leq 2n-2$,\\
&& $[e_i,e_{2n+1-i}]=(-1)^i e_{2n}$,& $2\leq i \leq n,$\\
 && $[e_i,x]=i e_i$, & $1\leq i \leq 2n-1$, & $[e_{2n},x]=(2n+1) e_{2n},$ \\ && $[e_i,y]=e_i$, & $1\leq i \leq 2n-1$, & $[e_{2n},y]=-[y,e_{2n}]=2 e_{2n}$.  \\
\hline
\end{longtable}

From the results of the works \cite{AbdAda} and \cite{KayKhud1}, we have the description of $\frac 1 2$-derivations of solvable Lie algebras with naturally graded filiform nilradicals $n_{n,1}$ and $Q_{2n}$ as follows:
\begin{longtable}{ll}
\hline
$\mathfrak{s}^{1}_{n,1}(2):$ & $\begin{array}{ll}
 \multicolumn{2}{l} {D(x)=\alpha_{1}x+\sum\limits_{i=2}^{n-1}(i-1)\alpha_{i+1}e_i+\delta e_n,} \\[1mm]
  D(e_1)=\sum\limits_{i=1}^{n}\alpha_ie_i, & \\[1mm]
  D(e_i)=\alpha_1e_i, & 2\leq i\leq n. \\[1mm]
 \end{array}$ \\
 \hline
 $\mathfrak{s}^{1}_{n,1}(\beta)_{\beta\neq 2}:$ & $\begin{array}{ll}
 \multicolumn{2}{l} {D(x)=\alpha_{1}x+\sum\limits_{i=2}^{n-1}(i-3+\beta)\alpha_{i+1}e_i+\delta e_n,}\\[1mm]
  D(e_1)=\alpha_1e_1+\sum\limits_{i=3}^{n}\alpha_ie_i, & \\[1mm]
  D(e_i)=\alpha_1e_i, &  2\leq i\leq n. \\[1mm]
 \end{array}$ \\ \hline
 $\mathfrak{s}^{2}_{n,1}:$ & $\begin{array}{ll}
  \multicolumn{2}{l} {D(x)=\alpha_{1}x+\sum\limits_{i=2}^{n-1}\alpha_{i+1}e_i+\delta e_n,} \\[1mm]
  D(e_1)=\alpha_1e_1+\sum\limits_{i=3}^{n}\alpha_ie_i, & \\[1mm]
  D(e_i)=\alpha_1e_i, &  2\leq i\leq n.\\[1mm]
  \end{array}$ \\ \hline
$\mathfrak{s}^{3}_{n,1}:$ & $\begin{array}{ll}
  \multicolumn{2}{l} {D(x)=\alpha_{1}x+\sum\limits_{i=3}^{n-1}(i-2)\alpha_{i+1}e_i+\delta e_n,} \\[1mm]
  D(e_1)=\alpha_1e_1+\sum\limits_{i=3}^{n}\alpha_ie_i, & \\[1mm]
  D(e_i)=\alpha_1e_i, & 2\leq i\leq n. \\[1mm]
  \end{array}$
  \\ \hline
  $\mathfrak{s}^{4}_{n,1}(\alpha_3, \dots, \alpha_{n-1}):$ & $\begin{array}{ll}
  \multicolumn{2}{l} {D(x)=a_1x+\sum\limits_{i=2}^{n-1}(a_{i+1}+\sum\limits_{t=3}^{i-1}\alpha_{t}a_{i-t+2})e_i+\delta e_n,}\\[1mm]
  D(e_1)=a_1e_1+\sum\limits_{i=3}^{n}a_ie_i, & \\[1mm]
  D(e_i)=a_1e_i, & 2\leq i\leq n. \\[1mm]
  \end{array}$ \\ \hline
$\mathfrak{s}_{n,2}:$ & $\begin{array}{ll}
D(x_1) = \alpha x_1 + (n-2)\beta e_n, \\
D(x_2) = \alpha x_2 + \beta e_n, \\
D(e_k) = \alpha e_k, & 1 \leq k  \leq n.\end{array}$
  \\ \hline
\end{longtable}

The description of $\frac 1 2$-derivations of solvable Lie algebras with naturally graded filiform nilradical $Q_{2n}$:
\begin{longtable}{ll}
\hline
$\tau_{2n,1}^1(\alpha)$ and $\tau_{2n,1}^2$: & $\begin{array}{ll}
 \multicolumn{2}{l} {D(x)=a x+(2n-3)b e_{2n-1}+ce_{2n},} \\[1mm]
  D(e_1)=ae_1, &  D(e_2)=a e_2+be_{2n}, \\[1mm]
 D(e_i)=ae_i, & 3\leq i \leq 2n. \\[1mm]
  
 \end{array}$ \\
 \hline
 $\tau_{2n,1}^3(\alpha_4, \alpha_6, \dots \alpha_{2n-2})$: & $\begin{array}{ll}
 D(x)=ax+ce_{2n},& \\[1mm]
  D(e_1)=a e_1, &  D(e_2)=a e_2+be_{2n},\\[1mm]
   D(e_i)=a e_i, & 3\leq i\leq 2n.\\[1mm]
 \end{array}$ \\ \hline
  $\tau_{2n,2}$: & $\begin{array}{ll}
 D(x_1)=ax_1+(2n+1)be_{2n},&\\[1mm]
 D(x_2)=ax_2+2be_{2n},& \\[1mm]
 D(e_i)=ae_i, & 1\leq i\leq 2n. \\[1mm]
 \end{array}$ \\ \hline
\end{longtable}

In the following theorem, we describe local $\frac 1 2$-derivations of solvable Lie algebras with naturally graded filiform nilradicals $n_{n,1}$ and $Q_{2n}.$

\begin{thrm}\label{locsolfil} Any local $\frac 1 2$-derivation
of the algebras $\mathfrak{s}^{1}_{n,1}(\beta),$  $\mathfrak{s}^{2}_{n,1},$ $\mathfrak{s}^{3}_{n,1},$ 
$\mathfrak{s}^{4}_{n,1}(\alpha_3, \dots, \alpha_{n-1}),$ $\mathfrak{s}_{n,2}$ $\tau_{2n,1}^1(\alpha),$  $\tau_{2n,1}^2$ $\tau_{2n,1}^3(\alpha_4,\alpha_6,\dots,\alpha_{2n-2})$ and 
$\tau_{2n,2}$ has the following form:
\begin{itemize}
\item for the algebra $ \mathfrak{s}^{1}_{n,1}(2):$
\begin{equation}\label{eq_11}\begin{array}{llll}             \Delta(x)=b_1x+\sum\limits_{j=2}^{n}b_{j}e_{j}, & \Delta(e_1)=b_1e_1+\sum\limits_{j=2}^nc_{j}e_j,              & \Delta(e_i)=d e_i, & 2\leq i\leq n.
              \end{array}\end{equation}
\item for the algebras $\mathfrak{s}^{1}_{n,1}(\beta)_{\beta\neq2},\ \mathfrak{s}^{2}_{n,1},\ \mathfrak{s}^{4}_{n,1}(\alpha_3, \dots, \alpha_{n-1}):$
\begin{equation*}\begin{array}{llll}
\Delta(x)=b_{1}x+\sum\limits_{j=2}^{n}b_{j}e_{j},&
        \Delta(e_1)=b_{1}e_1+\sum\limits_{j=3}^nc_{j}e_j, &
             \Delta(e_i)=d e_i, & 2\leq i\leq n.
\end{array}\end{equation*}

\item for the algebra $ \mathfrak{s}^{3}_{n,1}:$
\begin{equation*}\begin{array}{llll}
\Delta(x)=b_{1}x+\sum\limits_{j=3}^{n}b_{j}e_{j},& 
        \Delta(e_1)=b_{1}e_1+\sum\limits_{j=3}^nc_{j}e_j,&
 \Delta(e_i)=d e_i,& 2\leq i\leq n.
\end{array}\end{equation*}

\item for the algebra $ \mathfrak{s}_{n,2}:$
\begin{equation*}\begin{array}{llll}
               \Delta(x_1)=ax_1+(n-2)be_{n},&
              \Delta(x_2)=ax_2+be_{n},&
              \Delta(e_i)=ae_i, & 2\leq i\leq n.
\end{array}\end{equation*}

\item for the algebras $\tau_{2n,1}^1(\alpha),\ \tau_{2n,1}^2:$ \begin{equation*}\begin{array}{lll} 
\Delta(x)=ax+be_{2n-1}+ce_{2n}, & \Delta(e_1)=ae_1, \\
\Delta(e_2)=ae_2+de_{2n},&              \Delta(e_i)=ae_i,& 3\leq i \leq 2n.
              \end{array}\end{equation*}
\item for the algebras $\tau_{2n,1}^3(\alpha_4,\alpha_6,\dots,\alpha_{2n-2}):$
\begin{equation*}\begin{array}{lll}
\Delta(x)=ax+be_{2n}, & \Delta(e_1)=ae_1, & \\
    \Delta(e_2)=ae_2+ce_{2n}, & \Delta(e_i)=ae_i, & 3\leq i \leq 2n.
\end{array}\end{equation*}

\item for the algebra $ \tau_{2n,2}:$
\begin{equation*}\begin{array}{llll}
               \Delta(x_1)=ax_1+be_{2n},&
              \Delta(x_2)=ax_2+ce_{n}, &
              \Delta(e_i)=ae_i, & 1\leq i\leq 2n.
\end{array}\end{equation*}

\end{itemize}
\end{thrm}
\begin{proof} First, we give the proof of the theorem for the algebra $\mathfrak{s}^{1}_{n,1}(2).$ 
Let $\Delta$ be an arbitrary local $\frac 1 2$-derivation on $\mathfrak{s}^{1}_{n,1}(2).$
%and $\mathfrak{C}=(c_{i,j})$ be its matrix on the basis $\{x, e_1, e_2, \dots, e_n\}.$ 
%\[\mathfrak{C}=\left(\begin{array}{ccccccccc}
 %               c_{1,1}   & c_{1,2}        & \cdots & c_{1,n-1}   & c_{1,n}       \\
  %              c_{2,1}   & c_{2,2}        & \cdots & c_{2,n-1}   & c_{2,n}       \\
   %             c_{3,1}   & c_{3,2}        & \cdots & c_{3,n-1}   & c_{3,n}       \\
    %            \cdots    & \cdots         & \cdots & \cdots      & \cdots          \\
     %           c_{n-1,1} & c_{n-1,2}      & \cdots & c_{n-1,n-1} & c_{n-1,n}       \\
      %          c_{n,1}   & c_{n,2}        & \cdots & c_{n,n-1}   & c_{n,n}     \\
       %        \end{array}            \right)\]
Since $\Delta$ is a local $\frac 1 2$-derivation, then for any element $z=z_{1}x+\sum\limits_{i=1}^{n}z_{i+1}e_{i}\in \mathfrak{s}^{1}_{n,1}(2),$
there exists a $\frac{1}{2}$-derivation $D_z$, such that $\Delta(z)=D_z(z).$
%By Theorem \ref{halfderiv1}, the determinant $D_z$ has the following matrix form:
%\[\mathfrak{C}_z=\left(              \begin{array}{cccccccccc}
%\alpha_1^z                & 0              & 0            & 0           &  \cdots  &  %0        &0          & 0   \\
%  0                       & \alpha_1^z     & 0            & 0           &  \cdots  &  %0        &0          & 0 \\
%\alpha_3^z                &\alpha_2^z      & \alpha_{1}^z & 0           &  \cdots  &  %0        &0          & 0 \\
%2\alpha_4^z               &\alpha_3^z      & 0            & \alpha_1^z  &  \cdots  &  %0        &0          & 0 \\
% \cdots                   & \cdots         & \cdots       & \cdots      & \cdots   &  %\cdots   &\cdots     & \cdots \\
%(n-3)\alpha_{n-1}^z       & \alpha_{n-2}^z & 0            & 0           &  \cdots  %&\alpha_1^z &0          &  \\
%(n-2)\alpha_{n}^z         &\alpha_{n-1}^z  & 0            & 0           & \cdots   &  %0        &\alpha_1^z & 0   \\
%\delta_n^z                & \alpha_{n}^z   & 0            & 0           & \cdots   &  %0        &0          & \alpha_1^z\\   \end{array}            \right)\]
Choosing subsequently $z=x, z=e_1, \dots, z=e_n,$ we obtain the general form of $\Delta$ as follows:
\begin{equation*}\begin{array}{llll}             \Delta(x)=b_1x+\sum\limits_{j=2}^{n}b_{j}e_{j}, & \Delta(e_1)=\sum\limits_{j=1}^nc_{j}e_j,              & \Delta(e_i)=a_i e_i, & 2\leq i\leq n.
              \end{array}\end{equation*}

%the matrix $\mathfrak{C}=(c_{i,j})$ on basis $\{x, e_1, e_2, \dots, e_n\}:$ 
%\begin{equation}\label{eq1}
%\mathfrak{C}=\left(              \begin{array}{ccccccccc}                c_{1,1}   & 0              & 0        & 0         & \cdots   & 0          & 0         & 0 \\                0         & c_{2,2}        & 0        & 0         & \cdots   & 0          & 0         & 0 \\                c_{3,1}   & c_{3,2}        & c_{3,3}  & 0         & \cdots   & 0          & 0         &  0 \\                c_{4,1}   & c_{4,2}        & 0        & c_{4,4}   & \cdots   & 0          & 0         &  0 \\                \cdots    & \cdots         & \cdots   & \cdots    & \cdots   & \cdots     & \cdots    & \cdots \\                c_{n-1,1} &c_{n-1,2}       & 0        & 0         & \cdots   & c_{n-1,n-1}& 0         & 0\\                c_{n,1}   & c_{n,2}        & 0        & 0         & \cdots   & 0          &c_{n,n}    & 0\\                c_{n+1,1} & c_{n+1,2}      &  0       & 0         & \cdots   & 0          & 0         & c_{n+1,n+1} \\                   \end{array}                   \right)\end{equation}

Now consider $\Delta(e_2+e_i) = \Delta(e_2)+\Delta(e_i) = a_{3} e_2 + a_{i}e_i.$
On the other hand, there exist $\frac{1}{2}$-derivation $D_{e_2+e_i}$ such that $$\Delta(e_2+e_i) = D_{e_2+e_i}(e_2+e_i) = \alpha_{e_2+e_i}e_2+\alpha_{e_2+e_i}e_i.$$

Comparing  the coefficients of $e_2$ and $e_i$, we obtain that $a_{i}  =a_{3} $ for all $3 \leq i \leq n+1.$

Similarly, considering $$\Delta(x+e_1) = \Delta(x)+\Delta(e_1) = b_{1} x + c_{1}e_1 + \sum\limits_{j=2}^{n}(b_{j}+c_{j})e_{j},$$ 
and 
$$\Delta(x+e_1) = D_{x+e_1}(x+e_1) = \alpha_{x+e_1} x + \alpha_{x+e_1}e_1 + \sum\limits_{j=2}^{n-1}(\alpha_{j, x+e_1} + (j-1)\alpha_{j+1,x+e_1})e_{j} + (\alpha_{n, x+e_1} + \delta_{x+e_1})e_n,$$
we obtain that $c_1=b_1.$

Now, we show that any liner operator of the form 
\begin{equation*}\begin{array}{llll}             \Delta(x)=b_1x+\sum\limits_{j=2}^{n}b_{j}e_{j}, & \Delta(e_1)=b_1e_1+\sum\limits_{j=2}^nc_{j}e_j,              & \Delta(e_i)=d e_i, & 2\leq i\leq n,
              \end{array}\end{equation*}
is a local $\frac 1 2$-derivation on $\mathfrak{s}^{1}_{n,1}(2).$   
Considering $\Delta(z)=D_z(z)$ for any $z=z_{1}x+\sum\limits_{i=1}^{n}z_{i+1}e_{i},$
 we obtain the following  system of equality:
\begin{equation}\label{eq2}
  \begin{array}{lll}
    b_1z_1 &=\alpha_{1,z}z_1, \\
    b_1z_2 &=\alpha_{1,z}z_2, \\
    %b_{2}z_1+c_{2}z_2+az_{3} &=\alpha_3^zz_1+\alpha_{2}^zz_2+\alpha_{1}^zz_{3}, \\
    %b_{n-2}z_1+c_{n-2}z_{2}+a z_{n-1} &=(n-3)\alpha_{n-1}^zz_1+\alpha_{n-2}^zz_{2}+\alpha_{1}^zz_{n-1}, \\
    b_{i-1}z_1+c_{i-1}z_{2}+dz_{i} &=(i-2)\alpha_{i,z} z_1+\alpha_{i-1,z} z_{2}+\alpha_{1,z}z_{i}, & 3 \leq i \leq n,\\
    b_{n}z_1+c_{n}z_2+d z_{n+1}&=\delta_z z_1+\alpha_{n,z} z_2+\alpha_{1,z}z_{n+1}.
  \end{array}
\end{equation}

Let us consider the following cases:

\begin{itemize}\item[ Case 1.] Let $z_1\neq 0,$ then one can take 
\begin{equation*}\begin{split}
\alpha_{1,z}&=b_{1},\\
%\alpha_2^z&=0,\\
%\alpha_3^z&=\frac{c_{3,1}z_1+c_{3,2}z_2+c_{3,3}z_3-c_{1,1}z_3}{z_1},\\
\alpha_{i,z}&=\frac{b_{i-1}z_1+(c_{i-1}-\alpha_{i-1,z})z_2+(d-b_1)z_i}{(i-2)z_1},\quad 3\leq i\leq n,\\
\delta_z&=\frac{b_{n}z_1+(c_{n}-\alpha_{n,z})z_2+(d-b_1)z_{n+1}}{z_1}.
\end{split}\end{equation*}

\item[ Case 2.] Let $z_1=0$ and $z_2\neq0,$ then we may choose
\begin{equation*}\begin{split}
\alpha_{1,z}&=b_{1},\\
\alpha_{i-1,z}&=\frac{c_{i-1}z_2+(d-b_{1})z_i}{z_2},\quad 3\leq i\leq n+1.\\
\end{split}\end{equation*}

\item[ Case 3:] Let $z_1=z_2=\dots=z_{2}=0$, then we may take $\alpha_{1,z}=b_{1}.$ 
\end{itemize}

Thus, we have that for any $b_{1}, b_{2}, \dots, b_n, c_2, \dots, c_n, d$ and  $z_1, z_2, \dots, z_{n+1},$ there exist $\alpha_{i,z}$ and $\delta_z$, such that the equalities \eqref{eq2} hold. Hence any linear transformation of the form \eqref{eq_11} is a local $\frac 1 2$-derivation on $\mathfrak{s}^{1}_{n,1}(2).$

Now, we give the proof of the theorem for the algebra $\mathfrak{s}_{n,2}.$ Namely, we show that any local $\frac 1 2$-derivation on the algebra $\mathfrak{s}_{n,2}$ is a $\frac 1 2$-derivation.

Let  $\Delta$ be a local $\frac{1}{2}$-derivation on $\mathfrak{s}_{n,2},$ then 
considering the equalities $$\Delta(x_j)=D_{x_j}(x_j),\quad \Delta(e_i)=D_{e_i}(e_i),\quad 1\leq i\leq n,\quad 1\leq j\leq2,$$ we have
that
$$\begin{array}{ll}
\Delta(x_1)=d_{1}x_1+b_{1}e_n,&\Delta(x_2)=d_{2}x_2+b_{2}e_n,\\[2mm]
\Delta(e_i)=a_{i}e_i,&1 \leq i \leq n.
\end{array}$$

Now consider $\Delta(e_1+e_i)=a_{1}e_1+a_{i}e_i$ for $2\leq i\leq n.$ On the other hand,
\begin{equation*}
\Delta(e_1+e_i)=D_{e_1+e_i}(e_1+e_i)=D_{e_1+e_i}(e_1)+D_{e_1+e_i}(e_i) =\alpha_{e_1+e_i}e_1+\alpha_{e_1+e_i}e_i.
\end{equation*}

Comparing the coefficients at the basis elements $e_1$ and $e_i,$ we get $\alpha_{e_1+e_i}=a_{1},$
 $\alpha_{e_1+e_i}=a_{i},$
which implies $a_{i}=a_{1}$ for $2\leq i\leq n.$

Similarly, considering
\begin{equation*}
\begin{split}
\Delta(e_1+x_1)&=a_{1}e_1+d_{1}x_1+b_{1}e_n\\
 &=D_{e_1+x_1}(e_1+x_1)=D_{e_1+x_1}(e_1)+D_{e_1+x_1}(x_1)\\
                    & =\alpha_{e_1+x_1}e_1+\alpha_{e_1+x_1}x_1+(n-2)\beta_{e_1+x_1}e_n,\\
\end{split}
\end{equation*}
and
\begin{equation*}
\begin{split}
\Delta(e_1+x_2)&=a_{1}e_1+d_{2}x_2+b_{2}e_n\\
 &=D_{e_1+x_2}(e_1+x_2)=D_{e_1+x_2}(e_1)+D_{e_1+x_2}(x_2)\\
                    & =\alpha_{e_1+x_2}e_1+\alpha_{e_1+x_2}x_2+\beta_{e_1+x_2}e_n,\\
\end{split}
\end{equation*}
we get that 
$d_{1}=d_{2}=a_{1}.$

Finally, the equality
\begin{equation*}
\begin{split}
\Delta(x_1-(n-2)x_2)&=a_{1}x_1-(n-2)a_{1}x_2+(b_{1}-(n-2)b_{2})e_n\\
 &=D_{x_1-(n-2)x_2}(x_1-(n-2)x_2)=D_{x_1-(n-2)x_2}(x_1)-(n-2)D_{x_1-x_2}(x_2)\\
                    & =\alpha_{x_1-(n-2)x_2} x_1 -(n-2)\alpha_{x_1-(n-2)x_2}x_2,
\end{split}
\end{equation*}
implies
$b_{1}=(n-2)b_{2}.$

Thus, we obtain that the local $\frac{1}{2}$-derivation $\Delta$ has the following form:
$$\begin{array}{ll}
\Delta(x_1)=a_{1}x_1+(n-2)b_{1}e_n,&\Delta(x_2)=a_{1}x_2+b_{1}e_n,\\[2mm]
\Delta(e_i)=a_{1}e_i,&1 \leq i \leq n.
\end{array}$$

Hence, any local $\frac 1 2$-derivation on the algebra $\mathfrak{s}_{n,2}$ is a $\frac 1 2$-derivation.

The proofs are similar for the algebras $\mathfrak{s}^{1}_{n,1}(\beta)_{\beta\neq2},\ \mathfrak{s}^{2}_{n,1},\ \mathfrak{s}^{3}_{n,1},$ $\mathfrak{s}^{4}_{n,1}(\alpha_3, \alpha_4, \dots, \alpha_{n-1}),$ $\tau_{2n,1}^1(\alpha),$  $\tau_{2n,1}^2,$ $\tau_{2n,1}^3(\alpha_4,\alpha_6,\dots,\alpha_{2n-2})$ and 
$\tau_{2n,2}.$

\end{proof}

In the following table, we give the dimensions of the spaces of $\frac{1}{2}$-derivation and local $\frac{1}{2}$-derivations of solvable Lie algebras with naturally graded filiform nilradical:
$$\begin{tabular}{|p{4cm}|p{5cm}|p{5cm}|}
    \hline
  Algebra &  The dimensions of the space of $\frac{1}{2}$-derivations & The dimensions of the space of local $\frac{1}{2}$-derivations \\
  \hline
    $\mathfrak{s}^{1}_{n,1}(2)$  & $n+1$& $2n$ \\[1mm]
  \hline
    $\mathfrak{s}^{1}_{n,1}(\beta)_{\beta\neq2}$  & $n$& $2n-1$ \\[1mm]
  \hline
    $\mathfrak{s}^{2}_{n,1}$   & $n$& $2n-1$ \\[1mm]
  \hline
    $\mathfrak{s}^{3}_{n,1}$   & $n$& $2n-2$ \\[1mm]
  \hline
    $\mathfrak{s}^{4}_{n,1}(\alpha_3, \alpha_4, \dots, \alpha_{n-1})$   & $n$& $2n-1$ \\[1mm]
     \hline
    $\mathfrak{s}_{n,2}$   & $2$& $2$ \\[1mm]
  \hline
    $\tau_{2n,1}^1(\alpha),$  $\tau_{2n,1}^2$ & $3$ & $4$\\[1mm]
    \hline
    $\tau_{2n,1}^3(\alpha_4,\alpha_6,\dots,\alpha_{2n-2})$& $3$ & $3$\\[1mm]
    \hline
    $\tau_{2n,2}$& $2$ & $2$\\[1mm]
    \hline
  \end{tabular}$$

\begin{corollary}  Solvable Lie algebras $\mathfrak{s}^{1}_{n,1}(\beta),$  $\mathfrak{s}^{2}_{n,1},$ $\mathfrak{s}^{3}_{n,1},$ 
$\mathfrak{s}^{4}_{n,1}(\alpha_3, \dots, \alpha_{n-1}),$  $\tau_{2n,1}^1(\alpha),$ and $\tau_{2n,1}^2$ admit a local $\frac{1}{2}$-derivation which is not a $\frac{1}{2}$-derivation.
\end{corollary}

\begin{corollary}\label{Solloc}
Any local $\frac{1}{2}$-derivation on the algebras $\mathfrak{s}_{n,2},$  $\tau_{2n,1}^3(\alpha_4,\alpha_6,\dots,\alpha_{2n-2})$ and $\tau_{2n,2}$ is a $\frac{1}{2}$-derivation.
\end{corollary}

Now we investigate $2$-local $\frac{1}{2}$-derivations of the solvable Lie algebras with naturally graded filiform nilradicals and show that any $2$-local $\frac{1}{2}$-derivation of these algebras is a $\frac{1}{2}$-derivation.

\begin{thrm}\label{Sol2loc}
Any $2$-local $\frac{1}{2}$-derivation of the algebras  $\mathfrak{s}^{1}_{n,1}(\beta),$ $\mathfrak{s}^{2}_{n,1},$ $\mathfrak{s}^{3}_{n,1},$ $\mathfrak{s}^{4}_{n,1}(\alpha_3, \dots, \alpha_{n-1}),$
$\mathfrak{s}_{n,2},$ $\tau_{2n,1}^1(\alpha),$  $\tau_{2n,1}^2$ $\tau_{2n,1}^3(\alpha_4,\alpha_6,\dots,\alpha_{2n-2})$ and 
$\tau_{2n,2}$ is a $\frac{1}{2}$-derivation.
\end{thrm}

\begin{proof}  We will prove the theorem for the algebra \(\mathfrak{s}^{1}_{n,1}(2)\), and it is similar for the other previously given algebras.

First, we show that if $\nabla$ is a 2-local $\frac{1}{2}$-derivation
of \(\mathfrak{s}^{1}_{n,1}(2)\), such that $\nabla(x)=\nabla(e_{1})=0,$ then $\nabla\equiv0.$
Consider the pair of elements $x$ and $z=z_1x+\sum\limits_{i=1}^nz_{i+1}e_i.$ Then there exists a
$\frac{1}{2}$-derivation $D_{x,z},$ such that
$$\nabla(x)=D_{x,z}(x),\quad \nabla(z)=D_{x,z}(z).$$

Hence, from 
$$0=\nabla(x)=D_{x,z}(x)=\alpha_{1,x,z}x\sum\limits_{i=2}^{n-1}(i-1)\alpha_{i+1,x,z}e_i+\delta{x,z}e_n,$$
we have that  $\alpha_{1,x,z}=\alpha_{3,x,z}=\cdots=\alpha_{n,x,z}=\delta_{x,z}=0,$
which implies
\begin{equation}\label{2s1}
\nabla(z)=D_{x,z}(z)=z_2\alpha_{2,x,z}e_2.
\end{equation}

 Now, we consider the pair of elements $e_1$ and $z=z_1x+\sum\limits_{i=1}^nz_{i+1}e_i.$ Then there exists a
$\frac{1}{2}$-derivation $D_{e_1,z},$ such that
$$\nabla(e_1)=D_{e_1,z}(e_1),\quad \nabla(z)=D_{e_1,z}(z).$$

From 
$$0=\nabla(e_1)=D_{e_1,z}(e_1)=\sum\limits_{i=1}^{n}\alpha_{1,e_1,z}e_i,$$
we obtain  $\alpha_{1,e_1,z}=0$ for $1\leq i\leq n,$ which derives
\begin{equation}\label{2s2}
\nabla(z)=D_{e_1,z}(z)=z_1\delta_{e_1,z}e_n.
\end{equation}

From equalities \eqref{2s1} and \eqref{2s2}, we obtain that $\nabla(z)=0$ for any element $z \in \mathfrak{s}^{1}_{n,1}(2).$

Now let $\widetilde{\nabla}$ be any 2-local $\frac{1}{2}$-derivation of \(\mathfrak{s}^{1}_{n,1}(2)\). Take a $\frac{1}{2}$-derivation $D_{x,e_1}$ such that
\begin{equation*}
\widetilde{\nabla}(x)=D_{x,e_1}(x)\ \ \text{and} \ \ \widetilde{\nabla}(e_1)=D_{x,e_1}(e_1).
\end{equation*}
Set $\nabla_1=\widetilde{\nabla}-D_{x,e_1}.$ Then $\nabla_1$ is a 2-local
$\frac{1}{2}$-derivation such that $\nabla_1(x)=\nabla_1(e_1)=0,$ which implies   $\nabla_1\equiv0.$ Thus, $\widetilde{\nabla}=D_{x,e_1}$ is a $\frac{1}{2}$-derivation.
\end{proof}

%\begin{proof} Let $\nabla$ be a $2$-local $\frac{1}{2}$-derivation on $\mathfrak{s}_{n,2},$ such that $\nabla(x_1)=0.$
%Then for any element $\xi=\sum\limits_{i=1}^{n}\xi_ie_i+\xi_{n+1}x_1+\xi_{n+2}x_2\in \mathfrak{s}_{n,2},$ there exists a $\frac{1}{2}$-derivation $D_{x_1,\xi}(\xi)$, such that$$\nabla(x_{1})=D_{x_1,\xi}(x_1),\quad \nabla(\xi)=D_{x_1,\xi}(\xi).$$

%Hence,$$0=\nabla(x_{1})=D_{x_1,\xi}(x_1)=\alpha x_1+(n-2)\beta e_n,$$which implies,  $\alpha=\beta=0.$

%Consequently, from the description of the $\frac{1}{2}$-derivation $\mathfrak{s}_{n,2},$ we conclude that $D_{x_1,\xi}=0.$ Thus, we obtain that if $\nabla(x_1)=0,$ then $ \nabla\equiv0. $

%Let now $\nabla$ be an arbitrary $2$-local $\frac{1}{2}$-derivation of \(\mathfrak{s}_{n,2}\). Take a $\frac{1}{2}$-derivation $D_{x_1,\xi},$ such that
%\begin{equation*} \nabla(x_1)=D_{x_1,\xi}(x_1)\ \ \text{and} \ \ \nabla(\xi)=D_{x_1,\xi}(\xi).\end{equation*}

%Set $\nabla_1=\nabla-D_{x_1,\xi}.$ Then $\nabla_1$ is a $2$-local$\frac{1}{2}$-derivation, such that $\nabla_1(x_1)=0.$ Hence $\nabla_1(\xi)=0$ for all $\xi\in \mathfrak{s}_{n,2},$  which implies $\nabla=D_{x_1,\xi}.$Therefore, $\nabla$ is a $\frac{1}{2}$-derivation.\end{proof}

\subsection{Local and 2-local $\frac{1}{2}$-derivations on solvable Lie algebras with Heisenberg nilradical}
All solvable Lie algebras with Heisenberg nilradical were found in  \cite{rw93}.
In this subsection, we consider $(3n+2)$-dimensional solvable Lie algebra $\mathfrak L_{n,n+1}$ with the following multiplication table:
\begin{longtable}{lcllcllcll}
$[e_{n+i}, e_i]$&$=$&$e_{2n+1},$ &    $[e_i, x_{i}] $&$ =$&$ e_i,$  & $[e_{n+i}, x_{i}] $&$ =$&$ -e_{n+i},$  & $1 \leq i \leq n,$\\
& & & $[e_i, x_{n+1}] $&$=$&$ e_i,$  & $[e_{2n+1}, x_{n+1}] $&$= $&$e_{2n+1},$ & $1 \leq i \leq n.$
\end{longtable}

Note that nilradical of this algebra is the $(2n+1)$-dimensional
Heisenberg algebra $\mathfrak{h}_n$ with a basis $\{e_1, e_2, \dots, e_{2n+1}\}.$
In the following proposition, we give the description of $\frac 1 2$-derivations of the algebra $\mathfrak L_{n,n+1}.$

\begin{proposition} \cite{KayKhud1} \label{prop-Heisenberg}
Any $\frac 1 2$-derivation $D$ of the algebra $\mathfrak L_{n,n+1}$
has the form
 \begin{longtable}{lclllclllcl}
$D(e_k)$&$ =$&$ \alpha e_k, $&$ 1 \leq k  \leq 2n+1; $&\\[1mm]
$D(x_k) $&$= $&$\alpha x_k, $&$ 1 \leq k  \leq n; $&$
D(x_{n+1}) $&$= $&$\alpha x_{n+1} + \beta e_{2n+1}.$
\end{longtable}
\end{proposition}

\begin{thrm}\label{thm12}
 Any local $\frac{1}{2}$-derivation on the algebra $\mathfrak L_{n,n+1}$  is a $\frac{1}{2}$-derivation.
\end{thrm}

\begin{proof}

Let  $\Delta$ be a local $\frac{1}{2}$-derivation on $\mathfrak L_{n,n+1},$ then for any element $z\in \mathfrak L_{n,n+1},$ there exist a $\frac{1}{2}$-derivation $D_z$ such that 
$\Delta(z) = D_z(z).$ By Proposition \ref{prop-Heisenberg}, we obtain that $\frac{1}{2}$-derivation $D_z$ has the form 
\[\left\{\begin{array}{llllll}
D_z(e_i) = \alpha_{z} e_i, & 1 \leq i \leq 2n+1, \\
D_z(x_i) = \alpha_{z} x_i, & 1\leq i\leq n,& 
D_z(x_{n+1}) = \alpha_{z} x_{n+1} + \beta_{z} e_{2n+1}.
\end{array}\right.\]

Considering the equalities $$\Delta(x_j)=D_{x_j}(x_j),\quad \Delta(e_i)=D_{e_i}(e_i),\quad 1\leq i\leq 2n+1,\quad 1\leq j\leq n+1,$$ we conclude that
$$\begin{array}{lllll}
\Delta(e_i)=a_{i}e_i,&1 \leq i \leq 2n+1,\\[2mm]
\Delta(x_i)=b_{i}x_i,& 1 \leq i \leq n, &
\Delta(x_{n+1})=c x_{n+1}+d e_{2n+1}.
\end{array}$$

Now, consider $\Delta(e_1+e_i)=a_{1}e_1+a_{i}e_i$  for $2\leq i\leq 2n+1.$ 
On the other hand,
\begin{equation*}
\Delta(e_1+e_i)=D_{e_1+e_i}(e_1+e_i)=D_{e_1+e_i}(e_1)+D_{e_1+e_i}(e_i)=\alpha_{e_1+e_i}e_1+\alpha_{e_1+e_i}e_i.
\end{equation*}

Comparing the coefficients at the basis elements $e_1$ and $e_i,$ we get $\alpha_{e_1+e_i}=a_{1},$
 $\alpha_{e_1+e_i}=a_{i},$
which implies $a_{i}=a_{1}$ for $2\leq i\leq 2n+1.$

Similarly, considering 
\begin{equation*}
\begin{split}
\Delta(e_1+x_i)=a_{1}e_1+b_{i}x_i=D_{e_1+x_i}(e_1+x_i)=D_{e_1+x_i}(e_1)+D_{e_1+x_i}(x_i)=\alpha_{e_1+x_i}e_1+\alpha_{e_1+x_i}x_i,
\end{split}
\end{equation*}
for $2\leq i\leq n,$ and
\begin{equation*}
\begin{split}
\Delta(e_1+x_{n+1})&=a_{1}e_1+c x_{n+1}+d e_{2n+1}=D_{e_1+x_{n+1}}(e_1+x_{n+1})\\
 & = D_{e_1+x_{n+1}}(e_1)+D_{e_1+x_{n+1}}(x_{n+1}) =\alpha_{e_1+x_{n+1}}e_1+\alpha_{e_1+x_{n+1}} x_{n+1} + \beta_{e_1+x_{n+1}} e_{2n+1},\\
\end{split}
\end{equation*}
we derive $b_{i}=a_{1}$ for $1 \leq i \leq n$ and $c=a_1.$

Therefore, we obtain that any local $\frac{1}{2}$-derivation $\Delta$ has the following form:
$$\begin{array}{llll}
\Delta(e_i)=a_{1}e_i,&1 \leq i \leq 2n+1,\\[2mm]
\Delta(x_i)=a_{1}x_i,&1 \leq i \leq n,&
\Delta(x_{n+1})=a_{1}x_{n+1}+d e_{2n+1}.
\end{array}$$

Hence, every local $\frac{1}{2}$-derivation on $\mathfrak L_{n,n+1}$ is a $\frac{1}{2}$-derivation.
\end{proof}

\begin{thrm} \label{thm13}
Any $2$-local $\frac{1}{2}$-derivation of the algebra $\mathfrak L_{n,n+1}$ is a $\frac{1}{2}$-derivation.
\end{thrm}

\begin{proof} Let $\nabla$ be a $2$-local $\frac{1}{2}$-derivation on $\mathfrak L_{n,n+1}.$ Denote by $\mathcal{M} = \operatorname{span}\{e_1, \dots, e_{2n+1}, x_1, \dots, x_n\}.$ Since, any $\frac{1}{2}$-derivation of $\mathfrak L_{n,n+1}$ acts trivially on the subspace $\mathcal{M}$, that is $D(y) = \alpha y$ for all $y \in \mathcal{M},$ 
by Theorem \ref{thm1}, we conclude that $\nabla(y) = \alpha y$ for all $y \in \mathcal{M}.$

%Then, for the fixed element $x\in\mathcal{M}$ and for the arbitrary elements $y$ and $z$ of $\mathcal{M}$, we have that there exist  $\frac{1}{2}$-derivations $D_{x,y}$ and $D_{x,z}$, such that $$\nabla(x)=D_{x,y}(x) = \alpha_{x,y} x,\quad \nabla(y)=D_{x,y}(y) = \alpha_{x,y} y,$$ $$\nabla(x)=D_{x,z}(x) = \alpha_{x,z} x,\quad \nabla(z)=D_{x,z}(z) = \alpha_{x,z} z.$$

%Hence, we obtain $\alpha_{x,y} = \alpha_{x,z}.$ Due to the arbitrariness of $y$ and $z$, we get that $\nabla(y) = \alpha y$ for any element $y\in\mathcal{M}.$ 

Now, we consider any element $z \notin \mathcal{M}.$ Since $\nabla(\lambda x) = \lambda \nabla(x),$ then without loss of generality, we can suppose $z = x_{n+1} + m$ for some $m \in \mathcal{M}.$ For the arbitrary element $y\in \mathcal{M}$, we have that there exist  $\frac{1}{2}$-derivation $D_{y,z}$  such that
$$\alpha y = \nabla(y)=D_{y,z}(y) = \alpha_{y,z} y, \quad \nabla(z)=D_{y,z}(z) = \alpha_{y,z} x_{n+1}+\beta_{y,z}e_{2n+1} + \alpha_{y,z} m.$$

Then, we have that $\alpha_{y,z} = \alpha,$ which implies $\nabla(x_{n+1} + m) = \alpha x_{n+1}+\beta_m e_{2n+1} + \alpha m$ for any $m\in \mathcal{M}.$
Now, consider $\nabla(x_{n+1} + e_1) = \alpha x_{n+1}+\beta_{e_1}e_{2n+1} + \alpha e_1.$
On the other hand, from 
$$\begin{array}{lll}
    \nabla(z)&=D_{z,x_{n+1} + e_1}(z) &= \alpha x_{n+1}+\beta_{z,x_{n+1} + e_1}e_{2n+1} + \alpha m,\\
\nabla(x_{n+1} + e_1)&=D_{z,x_{n+1} + e_1}(x_{n+1} + e_1) &= \alpha x_{n+1}+\beta_{z,x_{n+1} + e_1}e_{2n+1} + \alpha e_1, \end{array}$$
we have that $\beta_m = \beta_{e_1}$ for any $m\in \mathcal{M}.$
Therefore, $\nabla$ is a $\frac{1}{2}$-derivation.
\end{proof}

%Let $\nabla$ be a $2$-local $\frac{1}{2}$-derivation on $L_{n,n+1},$ such that $\nabla(x_{n+1})=0.$Then for any element $\xi=\sum\limits_{i=1}^{2n+1}\xi_ie_i+\sum\limits_{i=1}^{n+1}\zeta_i x_i\in L_{n,n+1},$ there exists a $\frac{1}{2}$-derivation $D_{x_{n+1},\xi}(\xi)$, such that$$\nabla(x_{n+1})=D_{x_{n+1},\xi}(x_{n+1}),\quad \nabla(\xi)=D_{x_{n+1},\xi}(\xi).$$

%Hence,$$0=\nabla(x_{n+1})=D_{x_{n+1},\xi}(x_{n+1})=\alpha x_{n+1} + \beta e_{2n+1},$$which implies,  $\alpha=\beta=0.$

%Consequently, from the description of the $\frac{1}{2}$-derivation $L_{n,n+1},$ we conclude that $D_{x_{n+1},\xi}=0.$ Thus, we obtain that if $\nabla(x_{n+1})=0,$ then $ \nabla\equiv0. $

%Let now $\nabla$ be an arbitrary $2$-local $\frac{1}{2}$-derivation of \(L_{n,n+1}\). Take a $\frac{1}{2}$-derivation $D_{x_{n+1},\xi},$ such that \begin{equation*} \nabla(x_{n+1})=D_{x_{n+1},\xi}(x_{n+1})\ \ \text{and} \ \ \nabla(\xi)=D_{x_{n+1},\xi}(\xi). \end{equation*}

%Set $\nabla_1=\nabla-D_{x_{n+1},\xi}.$ Then $\nabla_1$ is a $2$-local $\frac{1}{2}$-derivation, such that $\nabla_1(x_{n+1})=0.$ Hence $\nabla_1(\xi)=0$ for all $\xi\in L_{n,n+1},$  which implies $\nabla=D_{x_{n+1},\xi}.$ Therefore, $\nabla$ is a $\frac{1}{2}$-derivation.

\subsection{Local and 2-local $\frac{1}{2}$-derivation on
solvable Lie algebras with abelian nilradical}

Now we consider solvable Lie algebras with abelian nilradical and maximal complementary vector space
(these algebras were found in \cite{nw94}).
It is known that the maximal dimension of complementary space for solvable Lie algebras with $n$-dimensional abelian nilpotent radical is equal to $n.$ Moreover, up to isomorphism there exists only one such solvable Lie algebra with the following multiplications:
$$\mathfrak{L}_n: \left[e_i, x_i\right]=e_i, \quad 1 \leq i \leq n,$$
where $\{e_1,   \dots, e_n, x_1,  \dots, x_n\}$ is a basis of $\mathfrak{L}_n.$

From the results of \cite{KayKhud1}, we have that 
 any $\frac 1 2$-derivation $D$ of the algebra $\mathfrak{L}_{n}$
has the form
\begin{equation} \label{eq5}
    D(e_i) = \alpha_i e_i,  \quad
D(x_i) = \alpha_i x_i + \beta_{i} e_i, \qquad 1 \leq i  \leq n.
\end{equation}

In the following theorem, we describe local 
$\frac{1}{2}$-derivations of the algebra $\mathfrak{L}_{n}.$

\begin{thrm}\label{locabel} Any local 
$\frac{1}{2}$-derivation of the algebra $\mathfrak{L}_{n}$
has the following form:
\begin{equation*}\begin{array}{lll}
    \Delta(e_i)&=a_{i}e_i, & 1\leq i\leq n,\\
    \Delta(x_i)&=b_{i}e_i+c_{i}x_i, & 1\leq i\leq n.
\end{array}
\end{equation*}
\end{thrm}
\begin{proof}
    The proof is similar to the proof of Theorem \ref{locsolfil}.
\end{proof}

From \eqref{eq5} and Theorem \ref{locabel}, it implies that $\dim(\mathrm{Der}_{\frac{1}{2}}(\mathfrak{L}_n)) =2n$ and $\dim(\mathrm{LocDer}_{\frac{1}{2}}(\mathfrak{L}_n) =3n.$ Thus, there exists a local $\frac{1}{2}$-derivation on $\mathfrak{L}_{n}$, which is not $\frac{1}{2}$-derivation.

%\begin{proof} Let us consider the linear operator $\Delta$ on $L_{n}$, such that
%$$\Delta\left(\sum\limits_{i=1}^n\xi_ie_i+\sum\limits_{i=1}^n\zeta_ix_i\right)=2\xi_1e_{1}+\zeta_1x_1.$$

%By Proposition \ref{abel}, it is not difficult to see that $\Delta$ is not a $\frac{1}{2}$-derivation.
%We show that, $\Delta$ is a local $\frac{1}{2}$-derivation on $L_{n}$.

%Consider the $\frac{1}{2}$-derivations $D_1$ and $D_2$ on the algebra $L_{n},$ defined as:
%\begin{equation*}\begin{split}
%D_1\left(\sum\limits_{i=1}^n\xi_ie_i+\sum\limits_{i=1}^n\zeta_ix_i\right)&=\xi_1e_{1}+\zeta_1x_1,\\
%%D_2\left(\sum\limits_{i=1}^n\xi_ie_i+\sum\limits_{i=1}^n\zeta_ix_i\right)&=\zeta_1x_1.
%\end{split}
%\end{equation*}

%Now, for any $\nu=\sum\limits_{i=1}^n\xi_ie_i+\sum\limits_{i=1}^n\zeta_ix_i,$ we shall find a $\frac{1}{2}$-derivation $D,$ such that $\Delta(\nu) = D(\nu).$

%If  $\xi_1=0,$ then
%$$\Delta(\nu)=0=D_2(\nu).$$

%If $\xi_1\neq 0,$ then setting
% $D=2D_1+tD_2,$ where $t=-\frac{\zeta_1}{\xi_1},$ we obtain that
%\begin{equation*}\begin{split}
 %\Delta(\nu)&=2\xi_1e_{1}+\zeta_1x_1=2\xi_1e_{1}+(2\zeta_1+t\xi_1)x_1=2(\xi_1e_{1}+\zeta_1x_1)+t\zeta_1x_1=\\
%                            &= 2D_1(\nu)+tD_2(\nu)=D(\nu).
%\end{split}\end{equation*}

%Hence, $\Delta$ is a local $\frac{1}{2}$-derivation.
%\end{proof}

\begin{thrm} \label{thm16}
Any $2$-local $\frac{1}{2}$-derivation of the algebra $\mathfrak{L}_{n}$ is a $\frac{1}{2}$-derivation.
\end{thrm}

\begin{proof}

Let $\nabla$ be a $2$-local $\frac{1}{2}$-derivation on $\mathfrak{L}_{n},$ such that $\nabla(q)=0,$ where $q=\sum\limits_{i=1}^{n}x_{i}.$
Then for any element $z=\sum\limits_{i=1}^{n}\xi_ie_i+\sum\limits_{i=1}^{n}\zeta_i x_i\in \mathfrak{L}_{n},$ there exists a $\frac{1}{2}$-derivation $D_{q,z}(z)$, such that
$\nabla(q)=D_{q,z}(q)$ and $ \nabla(z)=D_{q,z}(z).$
Then, from 
$$0=\nabla(q)=D_{q,z}(q)=\sum\limits_{i=1}^{n}\alpha_{i,q,z} x_i + \sum\limits_{i=1}^n\beta_{i,q,z} e_{i},$$
we obtain  $\alpha_{i,q,z}=\beta_{i,q,z}=0$ for $1\leq i\leq n,$
which implies $D_{q,z}=0.$
Thus, we obtain that if $\nabla(q)=0,$ then
$
\nabla\equiv0.
$

Now let $\widetilde{\nabla}$ be a any $2$-local $\frac{1}{2}$-derivation of \(\mathfrak{L}_{n}\).
Take a $\frac{1}{2}$-derivation $D_{q,z},$ such that
\begin{equation*}
\widetilde{\nabla}(q)=D_{q,z}(q)\ \ \text{and} \ \ \widetilde{\nabla}(z)=D_{q,z}(z).
\end{equation*}

Set $\nabla_1=\widetilde{\nabla}-D_{q,z}.$ Then $\nabla_1$ is a $2$-local
$\frac{1}{2}$-derivation, such that $\nabla_1(q)=0.$ Hence $\nabla_1(z)=0$ for all $z\in \mathfrak{L}_{n},$  which implies $\widetilde{\nabla}\equiv D_{q,z}.$
Therefore, $\widetilde{\nabla}$ is a
$\frac{1}{2}$-derivation.
\end{proof}

\subsection{Local and 2-local $\frac{1}{2}$-derivations on oscillator Lie algebras }\label{tpaoscl}

\begin{definition}
The oscillator Lie algebra $\mathfrak{L_{\lambda}}$ is a $(2n+2)$-dimensional Lie algebra with its canonical basis $\mathbb{B} = \{ e_{-1}, e_{0}, e_j, \check{e}_j \}_{j=1,\ldots,n}$
and the Lie brackets are given by
 $$[e_{-1},e_j] = \lambda_j\check{e}_j, \ \ \  [e_{-1},\check{e}_j] = -\lambda_je_j, \ \ \ [e_j,\check{e}_j] = e_0,$$
for $j = 1, \ldots, n$ and $\lambda = (\lambda_1, \ldots, \lambda_n) \in \mathbb{R}^{n}$ with $0 < \lambda_1 \leq \ldots \leq \lambda_n$.
\end{definition}

\begin{proposition}\cite{KayKhud1}\label{doscillator}
A linear map $D: \mathfrak{L}_{\lambda} \rightarrow \mathfrak{L}_{\lambda}$ is a $\frac{1}{2}$-derivation of the algebra $(\mathfrak{L}_{\lambda}, [-,-])$ if and only if
\begin{longtable}{lcl}
$D(e_{-1})$ &$ =$&$ \gamma e_{-1} + \mu e_0 - \sum\limits_{j=1}^{n} 2\lambda_j \alpha_{j} e_j - \sum\limits_{j=1}^{n} 2 \lambda_j \beta_{j} \check{e}_j$,\\
$D(e_0)$ &$ =$& $ \gamma e_0$,\\
$D(e_j)$ & $=$ & $\alpha_{j} e_0 + \gamma e_j, \ \ \ (j=1, \ldots, n)$,\\
$D(\check{e}_j)$ & $ =$ &$ \beta_{j} e_0 + \gamma \check{e}_{j}, \ \ \ (j=1, \ldots, n).$
\end{longtable}

\end{proposition}

\begin{thrm} \label{ldoscillator} Any local $\frac 1 2$-derivation on the oscillator Lie algebra $\mathfrak{L}_{\lambda}$  has the following form:
\begin{equation}\begin{array}{lll}\label{Lderosc}
              & \Delta(e_{-1})=\sum\limits_{j=-1}^{n}d_{j}e_{j}+\sum\limits_{j=1}^{n}d_{n+j}\check{e}_{j},\\
              & \Delta(e_0)=a_{0}e_0,\\
              & \Delta(e_i)=a_{i}e_0+b e_i,& 1\leq i\leq n,\\
              & \Delta(\check{e}_i)=c_{i}e_0+b \check{e}_i,& 1\leq i\leq n.\\
              \end{array}
\end{equation}
\end{thrm}
\begin{proof}
The proof is similar to the proof of Theorem \ref{locsolfil}.
\end{proof}

From Proposition \ref{doscillator} and Theorem \ref{ldoscillator}, it implies that $\dim(\frac{1}{2}\mathrm{Der}(\mathfrak{L_{\lambda}})) =2n+2$ and $\dim(\frac{1}{2}\mathrm{LocDer}(\mathfrak{L_{\lambda}})) =4n+4.$ Thus, there exists a local $\frac{1}{2}$-derivation on $\mathfrak{L_{\lambda}}$ which is not $\frac{1}{2}$-derivation.

\begin{thrm}\label{Sol2loc}
Any $2$-local $\frac{1}{2}$-derivation on the oscillator Lie algebra $\mathfrak{L}_{\lambda}$ is a $\frac{1}{2}$-derivation.
\end{thrm}

\begin{proof} The proof is similar to the proof of Theorem \ref{thm16}. The only difference here is that we consider a $2$-local $\frac{1}{2}$-derivation such that $\nabla(e_{-1})=0,$ and show that $\nabla(x)=0$ for any element $x.$ The other steps of the proof are the same. 
\end{proof}

\section{Local and 2-local $\frac{1}{2}$-derivations of some non-solvable Lie algebras}

In this section, we consider local and 2-local $\frac{1}{2}$-derivations for some finite-dimensional Lie algebras which is a semi-direct sum of semi-simple and solvable algebras, i.e., $\mathfrak{L} = \mathfrak{s} \ltimes \mathfrak{r}$. More precisely, we consider Lie algebras with three-dimensional simple parts and the solvable radical is an irreducible representation. Moreover, we consider local and 2-local $\frac{1}{2}$-derivations  Schr\"{o}dinger algebra $\mathcal{S}_{n}$ in $(n+1)$-dimensional space-time, which semi-simple part is a direct sum of three-dimensional algebra and orthogonal Lie algebra, solvable part is a  Heisenberg Lie
algebra.

\subsection{Local and 2-local $\frac{1}{2}$-derivation on semidirect products of ${\mathfrak {sl}}_2$ and irreducible modules.} 
It is known that if the  simple algebra $\mathfrak{sl}_2$ has an $(m+1)$-dimensional irreducible representation  $\mathfrak{r},$ where $m \geq 2$ (see, \cite{jac62}),
then we have the algebra $\mathfrak L^m = \mathfrak{sl}_2 \ltimes \mathfrak r$ with  the following multiplication:
{\small
\begin{longtable}{lcllcllcl}
$[e,f]$&$=$&$h,$ & $[h,e]$&$=$&$2e,$ & $[f,h]$&$=$&$2f,$\\
$[x_k,h]$&$=$&$(2k-m)x_k,$ & \multicolumn{6}{l}{$0 \leq k \leq m ,$}\\
$[x_k,f]$&$=$&$x_{k+1},$  & \multicolumn{6}{l}{$0 \leq k \leq m-1,$} \\
 $[x_k,e]$&$=$&$k(m+1-k)x_{k-1},$ & \multicolumn{6}{l}{$1 \leq k \leq m.$}
 \end{longtable}}

\begin{proposition}\cite{KayKhud1} \label{prop29}
Any $\frac 1 2$-derivation of $\mathfrak L^m (m\neq 2)$ is trivial and  $\frac 1 2$-derivation $D$ of $\mathfrak L^2$ has the following form
\begin{longtable}{lcllcllcl}
$D(e) $&$= $&$\alpha e -2\beta x_0,$ &
$D(f) $&$= $&$ \alpha f + \beta x_2,$ &
$D(h) $&$= $&$ \alpha h -2 \beta x_1,$ \\
$D(x_0) $&$= $&$\alpha x_0,$ & $D(x_1) $&$ = $&$ \alpha x_1,$ & $D(x_2) $&$= $&$\alpha x_2.$
\end{longtable}

 \end{proposition}

In the following theorems, we show that any local and 2-local $\frac{1}{2}$-derivation of the algebra $\mathfrak L^m$ is a $\frac{1}{2}$-derivation.  

\begin{thrm}
Any local $\frac{1}{2}$-derivation of the algebra $\mathfrak L^m$ is a $\frac{1}{2}$-derivation.
\end{thrm}
\begin{proof} According to Theorem \ref{thm1}, it is sufficient to prove the theorem for $m=2.$ Let  $\Delta$ be a local $\frac{1}{2}$-derivation on $\mathfrak L^2$, then for any element $z \in \mathfrak L^2$, there exist a $\frac{1}{2}$-derivation $D_z,$ such that $\Delta(z) = D_z(z).$
From Proposition \ref{prop29}, we have that  
\begin{longtable}{lcllcllcl}
$D_z(e) $&$= $&$\alpha_z e -2\beta_z x_0,$ &
$D_z(f) $&$= $&$ \alpha_z f + \beta_z x_2,$ &
$D_z(h) $&$= $&$ \alpha_z h -2 \beta_z x_1,$ \\
$D_z(x_0) $&$= $&$\alpha_{z} x_0,$ & $D_z(x_1) $&$ = $&$ \alpha_{z} x_1,$ & $D_z(x_2) $&$= $&$\alpha_{z} x_2.$
\end{longtable}

Considering the equalities $$\Delta(e)=D_{e}(e),\quad \Delta(f)=D_{f}(f),\quad \Delta(h)=D_{h}(h),\quad \Delta(x_j)=D_{x_j}(x_j),\quad 0\leq j\leq 2,$$ we obtain that 
$$\begin{array}{lll}
\Delta(e)=b_{1}e+c_1x_0, &
\Delta(f)=b_{2}f+c_2x_2, &
\Delta(h)=b_{3}h+c_3x_1, \\
\Delta(x_0)=b_{4}x_0, &
\Delta(x_1)=b_{5}x_1, &
\Delta(x_2)=b_{6}x_2. 
\end{array}
$$

Consider
$$\Delta(e+f+h)=b_{1}e+c_{1}x_0+b_{2}f+c_{2}x_2+b_{3}h+c_{3}x_1.$$

On the other hand,
\begin{equation*}
\begin{split}
\Delta(e+f+h)&=D_{e+f+h}(e+f+h)=D_{e+f+h}(e)+D_{e+f+h}(f)+D_{e+f+h}(h)\\
                    & =\alpha_{e+f+h} e -2\beta_{e+f+h} x_0+\alpha_{e+f+h} f +\beta_{e+f+h} x_2+\alpha_{e+f+h} h -2 \beta_{e+f+h} x_1.
\end{split}
\end{equation*}

Comparing the coefficients at the basis elements,  we get
$b_{1}=b_{2}=b_{3}$ and $c_{1}=c_{3}=-2c_{2}.$

Then, considering 
\begin{equation*}
\begin{split}
\Delta(e+x_1+x_2)&=b_{1}e-2c_2x_0+b_{5}x_1+b_{6}x_2\\
 &=D_{e+x_1+x_2}(e+x_1+x_2)=D_{e+x_1+x_2}(e)+D_{e+x_1+x_2}(x_1)+D_{e+x_1+x_2}(x_2)\\
                    & =\alpha_{e+x_1+x_2} e -2\beta_{e+x_1+x_2} x_0+\alpha_{e+x_1+x_2} x_1+\alpha_{e+x_1+x_2} x_2,\\
\end{split}
\end{equation*}
we obtain $b_{6}=b_{5}=b_{1}.$

Finally, from 
\begin{equation*}
\begin{split}
\Delta(f+x_0)&=b_{1}f+c_{2}x_2+b_{4}x_0\\
 &=D_{f+x_0}(f+x_0)=D_{f+x_0}(f)+D_{f+x_0}(x_{0})\\
                    & =\alpha_{f+x_0} f + \beta_{f+x_0} x_2+\alpha_{f+x_0} x_0,\\
\end{split}
\end{equation*}
we have that $b_{4}=b_{1}.$

Thus, we obtain that any local $\frac{1}{2}$-derivation $\Delta$ has the form:
$$\begin{array}{lll}
\Delta(e)=b_{1}e-2c_2x_0, &
\Delta(f)=b_{1}f+c_2x_2, &
\Delta(h)=b_{1}h-2c_2x_1, \\
\Delta(x_0)=b_{1}x_0, &
\Delta(x_1)=b_{1}x_1, &
\Delta(x_2)=b_{1}x_2. 
\end{array}
$$

Hence, every local $\frac{1}{2}$-derivation on $\mathfrak L^2$ is a $\frac{1}{2}$-derivation.
\end{proof}

\begin{thrm}
Any $2$-local $\frac{1}{2}$-derivation of the algebra $\mathfrak L^2$ is a $\frac{1}{2}$-derivation.
\end{thrm}

\begin{proof} Let $\nabla$ be a $2$-local $\frac{1}{2}$-derivation on $\mathfrak L^2,$ such that $\nabla(e)=0.$
Then for any element $z\in \mathfrak L^2,$ there exists a $\frac{1}{2}$-derivation $D_{e,z}$, such that
$\nabla(e)=D_{e,z}(e),$ $\nabla(z)=D_{e,z}(z).$
Then, from 
$$0=\nabla(e)=D_{e,z}(e)=\alpha_{e,z} e -2\beta_{e,z} x_0,$$
we get,  $\alpha_{e,z}=\beta_{e,z}=0,$ which implies $D_{e,z}=0.$
Thus, we obtain that if $\nabla(e)=0,$ then
$
\nabla\equiv0.
$

Let now $\nabla$ be an arbitrary $2$-local $\frac{1}{2}$-derivation of $\mathfrak L^2$.
Take a $\frac{1}{2}$-derivation $D_{e,z},$ such that
\begin{equation*}
\nabla(e)=D_{e,z}(e)\ \ \text{and} \ \ \nabla(z)=D_{e,z}(z).
\end{equation*}

Set $\nabla_1=\nabla-D_{e,z}.$ Then $\nabla_1$ is a $2$-local
$\frac{1}{2}$-derivation, such that $\nabla_1(e)=0.$ Hence $\nabla_1(z)=0$ for all $z\in \mathfrak L^2,$  which implies $\nabla=D_{e,z}.$
Therefore, $\nabla$ is a
$\frac{1}{2}$-derivation.
\end{proof}

%\subsection{Local and 2-local $\frac{1}{2}$ derivation on conformal %centrally extended Galilei algebras.}
%\begin{definition}\cite{klv22}  For every~$0< \ell \in \mathbb{N}-%\frac{1}{2}$, the
% conformal centrally extended Galilei algebra %$\widetilde{\mathfrak{g}}^{(\ell)}$  (it seems that it first appeared in %\cite{negro}) is generated by the following relations:
% \begin{longtable}{lcllcllcl}
%    $[h,e]$ & $=$ & $2e,$ & $[h,f]$ & $=$ & $-2f,$ & $[e,f]$ & $=$ & %$h,$\\
%    $[h,p_k]$ & $=$ & $2(\ell-k)p_k,$ & $[e,p_k]$ & $=$ & $kp_{k-1},$ & %$[f,p_k]$ & $=$ & $(2 \ell -k)p_{k+1},$ \\
%\multicolumn{9}{c}{$[p_k,p_{2 \ell - k}] =  (-1)^{k + \ell + \frac{1}%{2}}k! (2 \ell - k)! z,$}
%\end{longtable}
%  where  $k$  satisfies that  $0 \le k \le  2\ell$.
% \end{definition}
%\begin{proposition}\cite{klv22}
%Let $\widetilde{\mathfrak{g}}^{(\ell)}$ Galilei algebra.  Then any %$\frac 1 2$-derivation of $\widetilde{\mathfrak{g}}^{(\ell)}$ is trivial.
% \end{proposition}
%\begin{thrm}
%Any local and $2$-local $\frac{1}{2}$-derivation of Galilei algebra %$\widetilde{\mathfrak{g}}^{(\ell)}$ is a $\frac{1}{2}$-derivation.
%\end{thrm}
%\begin{proof}
%The proof is similar to the proof of Theorems \ref{Solloc} and %\ref{Sol2loc}.
%\end{proof}

\subsection{Local and 2-local $\frac{1}{2}$-derivation on Schr\"{o}dinger Lie algebras.} In this subsection, we describe local and 2-local $\frac{1}{2}$-derivation
 of the Schr\"{o}dinger algebra $\mathcal{S}_{n}$ in $(n+1)$-dimensional space-time.  

 The Schr{\"o}dinger algebra $\mathcal{S}_{n}$ is a finite-dimensional, non-semisimple and non-solvable Lie algebra, and it is the semidirect product Lie algebra
$$\mathcal{S}_{n}=(\mathfrak{sl}_{2}\oplus\mathfrak{so}_{n})\ltimes\mathfrak{h}_{n},$$
where $\mathfrak{sl}_{2}$ is the 3-dimensional simple Lie algebra, $\mathfrak{so}_{n}={\rm Span}_{\mathbb{C}}\{s_{kl}\mid1\leqslant k<l\leqslant n\}$ is the orthogonal Lie algebra and $\mathfrak{h}_{n}={\rm Span}_{\mathbb{C}}\{z,x_{i},y_{i}\mid1\leqslant i\leqslant n\}$ is the Heisenberg Lie algebra\cite{Dob1997}.

Thus, the Schr\"{o}dinger algebra $\mathcal{S}_{n}$ has a $\mathbb{C}$-basis
$$
\{e,f,h,z,x_{i},y_{i},s_{jk}(=-s_{kj})\mid1\leqslant i\leqslant n,1\leqslant j<k\leqslant n\},
$$
equipped with the following non-trivial commutation relations
$$\begin{array}{lll}
[e,f]=h, & [h,e]=2e,& [f,h]=2f,\\[1mm]
 [x_{i}, y_{i}]=z, & [h,x_{i}]=x_{i}, & [h, y_{i}]=-y_{i}, \\[1mm]
[e,y_{i}]=x_{i}, & [f,x_{i}]=y_{i},\\[1mm]
[s_{jk},x_{i}]=\delta_{ki}x_{j}-\delta_{ji}x_{k}, & [s_{jk},y_{i}]=\delta_{ki}y_{j}-\delta_{ji}y_{k},\\[1mm]
\multicolumn{2}{l}{[s_{jk},s_{lm}]=\delta_{lk}s_{jm}+\delta_{jm}s_{kl}+\delta_{mk}s_{lj}+\delta_{lj}s_{mk},}
\end{array}$$
where $\delta_{i,j}$ is the Kronecker Delta defined as $1$ for $i=j$ and as $0$ otherwise.

%\begin{remark}\label{rmk1}
  In particular, since the Schr{\"o}dinger algebra $\mathcal{S}_{1}$ does not have any term as $s_{jk}$ and 
in case of $n=2,$ we have the algebra $\mathcal{S}_{2}$ with a basis $\{f, h, e, z, x_1, y_1, x_2, y_2, s_{12}\},$ such that the multiplication is
$$\begin{array}{llll}
[e,f]=h, &  [h,e]=2e, & [f,h]=2f, \\[1mm]
[x_i, y_i] = z, &  [h, x_i] = x_i, &  [h, y_i] = -y_i, \\[1mm]
[e, y_i] = x_i, & [f, x_i] = y_i, \\[1mm]
[s_{12}, x_1] = -x_2, & [s_{12}, x_2] = x_1, & [s_{12}, y_1] = -y_2, & [s_{12}, y_2] = y_1.
\end{array}$$

In \cite{ytk} it is shown that there are no non-trivial $\frac 12$-derivations of the Schr\"{o}dinger algebra $\mathcal{S}_{n}$ for $n \neq 2$ and any 
$\frac 12$-derivation $D$ of $\mathcal{S}_{2}$ has the form 
\begin{equation}\label{dern=2}\begin{array}{llllll} D(e) = \alpha e, & D(f) = \alpha f, & D(h) = \alpha h, \\[1mm]
D(x_i) = \alpha x_i, & D(y_i) = \alpha y_i, &  1 \leq i \leq 2, \\[1mm]
 D(z) = \alpha z& D(s_{12}) = \alpha s_{12} + \beta z.\end{array}\end{equation}

\begin{thrm}
Any local and $2$-local $\frac{1}{2}$-derivation of Schr\"{o}dinger algebra $\mathcal{S}_{n}$ is a $\frac{1}{2}$-derivation.
\end{thrm}
\begin{proof}
The proof is similar to the proof of Theorems \ref{thm12} and \ref{thm13}.
\end{proof}

\begin{remark} Note that any local $\frac{1}{2}$-derivation of the algebra $\mathfrak L^2$ is a $\frac{1}{2}$-derivation, and $\mathfrak L^2$ does not have the non-trivial transposed Poisson structure \cite{KayKhud1}.

Any local $\frac{1}{2}$-derivation of the Schr{\"o}dinger algebra $\mathcal{S}_{2}$ is a $\frac{1}{2}$-derivation, but $\mathcal{S}_{2}$ 
has the non-trivial transposed Poisson structure \cite{ytk}.
\end{remark}

This remark motivates the following question:

{\bf Question.}  
Is there a Lie algebra that admits a local $\frac{1}{2}$-derivation, which is not a $\frac{1}{2}$-derivation, but does not have the non-trivial transposed Poison structure?


\begin{thebibliography}{99}


\bibitem{AbdAda} K.~Abdurasulov, J.~Adashev and S.~Eshmeteva, Transposed Poisson structure for solvable Lie algebra with filiform nilradical, (2024), arXiv:2401.04443.

\bibitem{AnCaGa3} J.~M.~Ancochea Berm\'{u}dez, R.~Campoamor-Stursberg,  L.~Garc\'{\i}a Vergnolle,  Solvable {L}ie algebras with naturally graded nilradicals and their invariants,  Journal of Physics A: Mathematical and Theoretical, 39 (2006), 6, 1339--1355.



%\bibitem{AlaYus2} A.~K.~Alauadinov, B.~B.~Yusupov,  Local derivations of conformal Galilei algebra, Communications in Algebra, (2024), https://doi.org/10.1080/00927872.2023.2301539.

\bibitem{AEK1} Sh.~Ayupov, A.~Elduque and K.~Kudaybergenov,   Local derivations and automorphisms of Cayley algebras, Journal of Pure and Applied Algebra,  227 (2023), 5, 107277.


\bibitem{AKR} Sh.~Ayupov, K.~Kudaybergenov, I.~Rakhimov,  2-Local derivations on finite-dimensional Lie algebras, Linear Algebra and its Applications, 474 (2015), 1--11.

\bibitem{Ayupov7} Sh.~Ayupov, K.~Kudaybergenov,
 Local derivations on finite-dimensional Lie algebras,  Linear Algebra and its Applications, 493 (2016), 381--398.

 
%\bibitem{AKA22} Sh.~A.~Ayupov,  K.~K.~Kudaybergenov,   A.~Allambergenov, Local and 2-local derivations on octonion algebras, Journal of algebra and its applications,  https://doi.org/10.1142/S0219498823501475.

\bibitem{AKO}  Sh.~A.~Ayupov,  K.~K.~Kudaybergenov, B.~A.~Omirov, Local and 2-local derivations and automorphisms on simple Leibniz algebras,  Bulletin of the Malaysian Mathematical Sciences Society, 43 (2020) 2199--2234.

%\bibitem{AKYu} Sh.~A.~Ayupov, K.~K.~Kudaybergenov, B.~B.~Yusupov, Local and 2-local derivations of $p$-filiform Leibniz algebras, Journal of Mathematical Sciences,  {\bf 245} (3),  (2020) 359--367.

\bibitem{AKYu1} Sh.~A.~Ayupov, K.~K.~Kudaybergenov, B.~B.~Yusupov,  Local and 2-Local Derivations of Locally Simple Lie Algebras, Journal of Mathematical Sciences,  278 (2024), 4, 613--622.

\bibitem{AKY} Sh.~Ayupov, K.~Kudaybergenov, B.~Yusupov,  2-Local derivations on generalized Witt algebras, Linear and Multilinear Algebra,  69 (2021), 16,  3130--3140.

\bibitem{AyuKhud} Sh.~Ayupov,  A.~Khudoyberdiyev, Local derivations on solvable Lie algebras, Linear and Multilinear Algebra,
69  (2021), 7,  1286--1301.

\bibitem{AyuKhudYus} Sh.~Ayupov,  A.~Khudoyberdiyev, B.~Yusupov, Local and 2-local derivations of solvable Leibniz algebras, International Journal of Algebra and Computation, 30 (2020), 6,  1185--1197.


\bibitem{AY1}  Sh.~Ayupov, B.~Yusupov,  2-local derivations of infinite-dimensional Lie algebras, Journal of Algebra and its
Applications,  19 (2020), 5, 2050100.


\bibitem{ChenXe} Q.-F.~Chen, Y.~He, 2-local derivations on the planar Galilean conformal algebra, International Journal of Mathematics, 34 (2023), 5, 2350023.



\bibitem{ChenWang} Z.~Chen, D.~Wang, 2-Local automorphisms of finite-dimensional simple Lie algebras,  Linear Algebra and its Applications,  486  (2015), 335--344.


\bibitem{CZZ} Y.~Chen, K.~Zhao, Y.~Zhao,   Local derivations on Witt algebras,
 Linear and Multilinear algebra, 70 (2022), 6, 1159--1172.

%\bibitem{Dora} E.~Dorado-Aguilar, R.~García-Delgado, E.~Martínez-Sigala, M.~C.~Rodríguez-Vallarte, G.~Salgado, Generalized derivations and some structure theorems for Lie algebras. J. Algebra Appl., {\bf 19}(2), (2020) 2050024.


\bibitem{Dob1997} V.~Dobrev, H.~Doebner, C.~Mrugalla, Lowest weight representations of the Schr\"{o}dinger algebra and generalized heat/Schr\"{o}dinger equations. Reports on Mathematical Physics 39, (1997), 2,  201--218.

%\bibitem{CP} C.~Duval, P.~Horvathy,  On Schrodinger superalgebras,  Journal of Mathematical Physics, {\bf 35}:5, (1994) 2516--2538.

%\bibitem{E11}  A.P.~Elisova, I.N.~Zotov, V.M.~Levchuk and G.S.~Suleymanova, Local automorphisms and local derivations of nilpotent matrix algebras, Izv. Irkutsk Gos. Univ., {\bf 4}:1   (2011) 9-19.
	
%\bibitem{BIK} B.~Ferreira, I.~Kaygorodov, K.~Kudaybergenov,  Local and 2-local derivations of simple n-ary algebras, Ricerche di Matematica, (2021).

\bibitem{BKL} B.~Ferreira, I.~Kaygorodov, V.~Lopatkin,
$\frac{1}{2}$-derivations of Lie algebras and transposed Poisson algebras, 
Revista de la Real Academia de Ciencias Exactas, Físicas y Naturales. Serie A. Matemáticas, 115 (2021), 3, 142. 

 \bibitem{fil1} V.Filippov,  On $\delta$-derivations of Lie algebras. Siberian Mathematical Journal, 39 (1998), 6, 1218--1230.
 
 
\bibitem{fil2} V.Filippov, $\delta$-derivations of prime Lie algebras. Siberian Mathematical Journal, 40 (1999), 1, 174--184.
 
 \bibitem{fil3} V.Filippov, On $\delta$-derivations of prime alternative and Mal’tsev algebras. Algebra and Logic, 39, (2000), 5, 354--358.
 
\bibitem{jac62} N.~Jacobson, Lie algebras, 
Interscience Tracts in Pure and Applied Mathematics, No. 10. Interscience Publishers (a division of John Wiley $\&$ Sons, Inc.), New York-London, 1962. ix+331 pp.

\bibitem{QX} Q.~Jiang, X.~Tang, 2-Local derivations on the Schrödinger-Virasoro algebra, Linear and Multilinear Algebra, (2023), DOI: 10.1080/03081087.2023.2176421.

%\bibitem{Jon}	B.E.~Johnson,  Local derivations on  $C^\ast$-algebras are derivations,  Transactions of the American Mathematical Society, 353, (2001)  313--325.


\bibitem{Kadison} R.V.~Kadison,   Local derivations, Journal of Algebra,  130, (1990), 494--509.

%\bibitem{Kac} V.~Kac, A.~Raina, Bombay lectures on highest weight representations of infinite-dimensional Lie algebras, World Sci.Singapore, 1987.


%\bibitem{Kay1} I.~Kaygorodov, $\delta$-superderivations of simple finite-dimensional Jordan and Lie superalgebras. Algebra and Logic,  49 (2010), 2, 130--144.

\bibitem{Kay2} I.~Kaygorodov, $\delta$-derivations of classical Lie superalgebras, Siberian Mathematical Journal, 50 (2009), 3, 434--449.



\bibitem{Kay3} I.~Kaygorodov, On $\delta$-derivations of simple finite-dimensional Jordan superalgebras, Algebra and Logic, 46 (2007), 5, 318--329.

\bibitem{Kay-n-ary} I.~Kaygorodov, $(n + 1)$-Ary derivations of semisimple Filippov algebras. Mathematical Notes, 96 (2014), 2, 208--216.

\bibitem{KayKhudYul} I.~Kaygorodov, K.~Kudaybergenov, I.~Yuldashev,
Local derivations of semi-simple Leibniz algebras. 
Communication in Mathematics. 30 (2022), 2, 1--12.


\bibitem{KayKhud1}  I.~Kaygorodov, A.~Khudoyberdiyev, Transposed Poisson structures on solvable and perfect Lie algebras, Journal of Physics A: Mathematical and Theoretical, 57 (2024), 3, 035205.

\bibitem{kk23}
I.~Kaygorodov, M.~Khrypchenko,
Transposed Poisson structures on Witt-type algebras,
 Linear Algebra and its Applications,  665  (2023),  196--210.	


\bibitem{kkg23}
I.~Kaygorodov, M.~Khrypchenko,
Transposed Poisson structures on generalized Witt algebras and Block Lie algebras,    Results in Mathematics, 78 (2023), 5,  186.

%\bibitem{KK7} Kaygorodov I.,  Khrypchenko M.,Transposed Poisson structures on the Lie algebra of upper triangular matrices, arXiv:2305.00727.

%\bibitem{kkinc}Kaygorodov I.,  Khrypchenko M.,Transposed Poisson structures on   Lie incidence algebras, arXiv:2309.00332. 



\bibitem{klv22} I.~Kaygorodov, V.~Lopatkin, Z.~Zhang, Transposed Poisson structures on Galilean and solvable Lie algebras,  Journal of Geometry and Physics, 187 (2023), 104781.

\bibitem{KayOpa} I.~Kaygorodov, E.~Okhapkina, $\delta$-derivations of semi-simple finite-dimensional structurable algebras, Journal of Algebra and its Applications, 13 (2014), 4, 1350130.

 
\bibitem{KOK22} K.~Kudaybergenov,  B.~Omirov, T.~Kurbanbaev,  Local derivations on solvable Lie algebras of maximal rank, Communications in Algebra,   50 (2022), 9, 1-11.

%\bibitem{KKY} K.K.~Kudaybergenov,  I.~Kaygorodov, I.~Yuldashev, Local derivations of semisimple Leibniz algebras, Communications in Mathematics {\bf 30}:2 (2022)  1–12.

%\bibitem{PatZas} Patera J., Zassenhaus H., The Higher Rank Virasoro Algebras. Commun. Math. Phys., (136): 1, (1991) 1-14.

\bibitem{Larson} D.~R.~Larson, A.~R.~Sourour,  Local derivations and local automorphisms of \(B(X)\),  Proceedings of Symposia in Pure Mathematics, 51 (1990), 187--194.


%\bibitem{Leger} G.~Leger, E.~Luks, Generalized derivations of Lie algebras, Journal of Algebra, 228 (2000), 1, 165--203.


%\bibitem{Liu} G.~Liu, Y.~Li, K.~Wang, Irreducible weight modules over the Schr\"{o}dinger Lie algebra in $(n+1)$ dimensional space-time, J. Algebra,  Vol. 575, (2021), 1--13.

\bibitem{nw94} J.~Ndogmo,  P.~Winternitz,
Solvable Lie algebras with abelian nilradicals, Journal of Physics A: Mathematical and Theoretical,   27,  (1994), 2, 405--423.

 \bibitem{rw93} J.~Rubin,   P.~Winternitz,
Solvable Lie algebras with Heisenberg ideals,
Journal of Physics A: Mathematical and Theoretical,
  26 (1993), 5, 1123--1138.

 \bibitem{Sem} \v{S}emrl P.,  Local automorphisms and derivations on $B(H),$ Proceedings of the American Mathematical Society, 125 (1997), 2677--2680.

\bibitem{sw05} L.~Šnobl,   P.~Winternitz,
A class of solvable Lie algebras and their Casimir invariants,
Journal of Physics A: Mathematical and Theoretical,
 38 (2005), 12, 2687--2700.


%\bibitem{Tang} X.~Tang, 2-Local derivations on the $W$-algebra $W(2,2)$, Journal of Algebra and Its ApplicationsVol. 20, No. 12, (2021) 2150237.

\bibitem{XMP} X.~Tang, M.~Xiao and P.~Wang, Local properties of Virasoro-like algebra, Journal of Geometry and Physics, 186 (2023), 104772.

\bibitem{Ver}  M.~Vergne,  Cohomologie des alg\`ebres de {L}ie nilpotentes. {A}pplication \`a  l\'etude de la vari\'et\'e des alg\`ebres de {L}ie nilpotentes, Bulletin de la Soci\'{e}t\'{e} Math\'{e}matique de France, 98 (1970), 81--116,
	
\bibitem{Yao} 	Y. F. Yao, 	 Local derivations on the Witt algebra in prime characteristic,  Linear and multilinear algebra, 70 (2022), 15, 2919--2933.

\bibitem{ytk} Ya.~Yang,  X.~Tang,  A.~Khudoyberdiyev, Transposed Poisson structures on Schrodinger algebra in $(n+1)$-dimensional space-time, (2023), arXiv:2303.08180.

\bibitem{YK} Y.~Yao, K.~Zhao, Local properties of Jacobson-Witt algebras,
Journal of Algebra, 586 (2021),  1110--1121.

	\bibitem{YC20}	Y.~Yu, Zh.~Chen,  Local derivations on Borel subalgebras of finite-dimensional simple Lie algebras,  Communications in Algebra, 48 (2020), 1,  1--10.



%\bibitem{YYZ07}		Y.~Yao, Y.~Ye, P.~Zhang,		Quiver Poisson algebras, Journal of  Algebra, 312 (2007), 2, 570--589.


    \bibitem{yh21}
L.~Yuan,  Q.~Hua,
$\frac{1}{2}$-(bi)derivations and transposed Poisson algebra structures on Lie algebras,  Linear and Multilinear Algebra, 70 (2022),   22, 7672--7701.


\bibitem{Zus}  P.~Zusmanovich,
  On $\delta$-derivations of Lie algebras and superalgebras, Journal of Algebra, 324,  (2010), 12, 3470--3486.
		




\end{thebibliography}
\end{document}